\documentclass[11pt,a4paper]{article}
\usepackage{amssymb,amsmath}
\usepackage{graphicx}
\newcommand{\trp}{{\sf\scriptsize T}}
\newcommand\eop{\hfill$\Box$\\}
\newtheorem{theorem}{Theorem}[section]
\newtheorem{lemma}{Lemma}[section]
\newtheorem{corollary}{Corollary}[section]
\newcommand\asto{\,\stackrel{\rm\scriptsize a.s.}{\longrightarrow}\,}
\newcommand\dto{\,\stackrel{\rm\scriptsize d}{\longrightarrow}\,}
\newcommand\Ifkt {\hskip.1em 1\hskip-.6em 1\hskip.1em}
\newcommand\thb{\overline{\theta}}
\newcommand\LSk{\widehat{\theta}_k^{\rm\scriptsize(LS)}}
\newcommand\MLk{\widehat{\theta}_k^{\rm\scriptsize(ML)}}
\newcommand\Ms{M_{{\rm\scriptsize s}*}}

\def\BM#1{{\mathchoice%
    {\mbox{\boldmath$#1$}}%
    {\mbox{\boldmath$#1$}}%
    {\mbox{\boldmath$\scriptstyle#1$}}%
    {\mbox{\boldmath$\scriptscriptstyle#1$}}}}
\begin{document}
\title{A $p$-step-ahead sequential adaptive algorithm\\ for D-optimal nonlinear regression design}
\author{Fritjof Freise$^1$, Norbert Gaffke$^2$ \ and \  Rainer Schwabe$^2$ \\[2ex]
\parbox{13cm}{\small $^1$University of Veterinary Medicine Hannover \ and \ 
$^2$University of Magdeburg}
}

\maketitle

\begin{abstract}
Under a nonlinear regression model with univariate response an algorithm 
for the generation of sequential adaptive designs is studied. 
At each stage, the current design is augmented by adding $p$ design points 
where $p$ is the dimension of the parameter of the model. 
The augmenting $p$ points are such that, at the current parameter estimate, 
they constitute the locally D-optimal design within the set of all saturated designs.
Two relevant subclasses of nonlinear regression models are focused on, which were considered in previous work
of the authors on the adaptive Wynn algorithm: firstly, regression
models satisfying the `saturated identifiability condition' and, secondly, generalized linear models.
Adaptive least squares estimators and adaptive maximum likelihood
estimators in the algorithm are shown to be strongly consistent and asymptotically normal,
under appropriate assumptions. 
For both model classes, if a condition of `saturated D-optimality' is satisfied,
the almost sure asymptotic D-optimality of the generated design sequence is implied by the strong consistency of 
the adaptive estimators employed by the algorithm. 
The condition states that there is a saturated design which is locally D-optimal at  the true parameter point
(in the class of all designs). 
\end{abstract} 

\section{Introduction}
\setcounter{equation}{0}
Sequential adaptive design and estimation in nonlinear regression models were considered by 
Lai and Wei \cite{Lai-Wei}, Lai \cite{Lai}, and Chen, Hu and Ying \cite{Chen-Hu-Ying}. 
In those fundamental contributions fairly general 
conditions on the adaptive design ensure consistency and asymptotic normality of adaptive least squares 
or maximum quasi-likelihood estimators. However, it remains open whether particular sequential adaptive design
schemes are covered, like the adaptive version of the algorithm of Wynn \cite{Wynn} for D-optimal design,
which we have called the `adaptive Wynn algorithm'. Pronzato \cite{Pronzato} was the first who studied the asymptotics of
the adaptive Wynn algorithm, that is, the asymptotic properties of the adaptive designs and adaptive least squares 
and maximum likelihood estimators under the algorithm. Crucial assumptions in that paper are a finite experimental region
and a condition of ``saturated identifiability'' (see below) on the regression model. 
Extensions of results in \cite{Pronzato} to any compact experimental region, and further results on the adaptive Wynn algorithm 
have been obtained by the authors in \cite{FF-NG-RS-18} and \cite{FF-NG-RS-19}. 
In the present paper a sequential adaptive design
algorithm is proposed and studied which we call ``$p$-step-ahead algorithm'' since at each step a batch of $p$ further design 
points is collected. For a special model a related concept of ``batch sequential design'' was employed by   
M\"uller and P\"otscher \cite{Mueller-Poetscher}. An idea of the algorithm was sketched 
by Ford, Torsney and Wu \cite{Ford-Torsney-Wu}, p.~570, in the introduction of their paper.    
Note that the adaptive Wynn algorithm collects one design point at each step and   
was therefore called ``one-step ahead algorithm'' in \cite{Pronzato}. 
Actually, in dimension $p=1$ both algorithms coincide.
When $p\ge2$, a practical advantage of the adaptive $p$-step-ahead algorithm over a strictly sequential $1$-step-ahead
sampling scheme like the adaptive Wynn algorithm might be that it allows some parallel response sampling
(batches of size $p$) and thus reduces the total duration of data collection. 

The paper is organized as follows. In Section 2 the general framework is outlined, and various conditions on the 
nonlinear regression model are introduced which will later be assumed for some results but not throughout. Some examples
of frequently used nonlinear models are discussed. In Section 3 the $p$-step-ahead algorithm is described.
In Section 4 some basic asymptotic properties of the design sequence generated by the algorithm are derived.
Sections 5 and 6 address consistency and asymptotic normality of the adaptive least squares and 
maximum likelihood estimators in the algorithm. An appendix contains supplementary results to two examples (parts A.1 and A.2 of the appendix)
and the proofs of the lemmas and theorems (parts A.3 and A.4). 

\section{General framework}
Let a nonlinear regression model be given with
univariate mean response function $\mu(x,\theta)$, $x\in{\cal X}$, $\theta\in\Theta$,
where ${\cal X}$ and $\Theta$ are the experimental region and the parameter space, respectively.
Also, a family of $\mathbb{R}^p$-valued functions $f_\theta$, $\theta\in\Theta$, defined on ${\cal X}$
is given such that the $p\times p$ matrix $f_\theta(x)\,f_\theta^\trp(x)$ constitutes the elemental 
information matrix of $x\in{\cal X}$ at $\theta\in\Theta$. Note that a vector $a\in\mathbb{R}^p$ is
written as a column vector and $a^\trp$ denotes its transposed which is a $p$-dimensional row vector.
An approximate design, for short: design, is a probability measure $\xi$ on ${\cal X}$ with finite support. 
The support of a design $\xi$ is denoted by ${\rm supp}(\xi)$, which is a nonempty finite subset of ${\cal X}$.
The weights $\xi(x)$ for $x\in{\rm supp}(\xi)$ are positive real numbers with 
$\sum_{x\in {\scriptsize\rm supp}(\xi)}\xi(x)\,=1$.
The information matrix of a design $\xi$ at $\theta\in\Theta$ is defined by
\begin{equation}
M(\xi,\theta)\,=\,\sum_{x\in {\scriptsize\rm supp}(\xi)}\xi(x)\,f_\theta(x)\,f_\theta^\trp(x),
\label{eq1-1}
\end{equation} 
which is a nonnegative definite $p\times p$ matrix. 
Throughout, as in \cite{FF-NG-RS-19}  the following basic conditions (b1) to (b4) are assumed.
\begin{itemize}
\item[(b1)]
The experimental region ${\cal X}$ is a compact metric space.
\item[(b2)] 
The parameter space $\Theta$ is a compact metric space.
\item[(b3)] 
The real-valued mean response function \ $(x,\theta)\mapsto \mu(x,\theta)$, defined on the Cartesian product space  
${\cal X}\times \Theta$, is continuous.
\item[(b4)]
The family $f_\theta$, $\theta\in\Theta$, of $\mathbb{R}^p$-valued functions on ${\cal X}$ satisfies:\\
(i) for each $\theta\in\Theta$ the image  $f_\theta({\cal X})$ \ spans $\mathbb{R}^p$;\\
(ii) the function \ $(x,\theta)\mapsto f_\theta(x)$, defined on   
${\cal X}\times \Theta$, is continuous.
\end{itemize}

More specific conditions will be employed later which, however, will not be assumed throughout. 
Next we formulate some  of them: condition (SI)
on ``saturated identifiabiliy'' as in \cite{FF-NG-RS-19}, condition (GLM) taking up
particular features of a generalized linear model
as in \cite{FF-NG-RS-18}, and a slightly stronger condition (${\rm GLM}^*$).

\vspace*{2ex}\noindent
{\em Condition {\rm(SI)}}\\   
For all pairwise distinct points $x_1,\ldots,x_p\in{\cal X}$ the $\mathbb{R}^p$-valued function
on $\Theta$,\\ 
$\theta\mapsto\bigl(\mu(x_1,\theta),\ldots,\mu(x_p,\theta)\bigr)^\trp$,
is an injection, that is, if $\theta,\theta'\in\Theta$ and $\mu(x_j,\theta)=\mu(x_j,\theta')$ for $j=1,\ldots,p$, 
then $\theta=\theta'$.\\

\noindent
{\em Condition {\rm (GLM)}}\\ 
$f_\theta(x)\,=\,\psi(x,\theta)\,f(x)$ \ for all $(x,\theta)\in {\cal X}\times \Theta$, where
$\psi\,:\,{\cal X}\times\Theta\longrightarrow(\,0\,,\,\infty)$ \ and \ $f:\,{\cal X}\longrightarrow\mathbb{R}^p$ \ 
are given continuous functions.\\

\noindent
{\em Condition $({\rm GLM}^*)$}\\ 
(GLM) holds and, moreover: $\Theta\subseteq\mathbb{R}^p$, $\mu(x,\theta)=G\bigl(f^\trp(x)\,\theta\bigr)$ for all 
$(x,\theta)\in{\cal X}\times\Theta$, where $G\,:\,I\longrightarrow\mathbb{R}$ is a continuously differentiable function
on an open interval $I\subseteq\mathbb{R}$ with positive derivative $G'>0$, and
$f^\trp(x)\,\theta\in I$ for all $(x,\theta)\in{\cal X}\times\Theta$.\\

A further condition refers to some given parameter point $\thb\in\Theta$, which
might be called a condition of ``saturated local D-optimality at $\thb$'', abbreviated by 
(SD$^*$)$(\thb)$. Note that a design $\xi^*_{\thb}$ is called locally D-optimal at $\thb$
if $\xi^*_{\thb}$ maximizes $\det M(\xi,\thb)$ over the set of all designs $\xi$.
A design is called saturated if its support size is equal to $p$.

\vspace*{2ex}\noindent
{\em Condition $({\rm SD}^*)(\thb)$}\\ 
There exists a locally D-optimal design at $\thb$ which is saturated.\\

For some results a weaker condition (SD)$(\thb)$ will be employed, which addresses the saturated designs 
maximizing the D-criterion locally at $\thb$ over the set of all saturated designs. 
For short, we call such designs ``locally D-optimal saturated designs at $\thb$''.
Note that a locally D-optimal saturated design at $\thb$ has uniform weights,   
since for any saturated design $\eta$ with support points $x_1,\ldots,x_p\in{\cal X}$
one gets from  (\ref{eq1-1}) 
\begin{equation}
\det M(\eta,\thb)\,=\,\Bigl(\prod_{j=1}^p\eta(x_j)\Bigr)\,
\bigl(\det\bigl[f_{\thb}(x_1),\ldots,f_{\thb}(x_p)\bigr]\bigr)^2, \label{eq1-2}
\end{equation}
and the product of the weights is maximized iff $\eta(x_j)=1/p$ for all $j=1,\ldots,p$. Thus, a locally D-optimal 
saturated design at $\thb$ is an equally weighted design on $p$ points $x_1^*,\ldots,x_p^*\in{\cal X}$ which
maximize 
$\bigl(\det\bigl[f_{\thb}(x_1),\ldots,f_{\thb}(x_p)\bigr]\bigr)^2$ over $x_1,\ldots,x_p\in{\cal X}$.
This motivates the iteration rule of the $p$-step-ahead algorithm, see (\ref{eq2-1}) in Section 3.\\  

\noindent
{\em Condition $({\rm SD})(\thb)$}\\ 
The information matrices $M(\eta^*,\thb)$ of all locally D-optimal saturated designs $\eta^*$ 
at $\thb$ coincide, and are thus equal to one matrix $\Ms(\thb)$, say.\\

As it is well-known, the locally D-optimal information matrix at $\thb$ is unique,
$M_*(\thb)$ say. Therefore, if condition (SD$^*$)$(\thb)$ holds then condition (SD)$(\thb)$
holds as well and $\Ms(\thb)=M_*(\thb)$.
There are several relevant nonlinear   
models which satisfy condition $({\rm SD}^*)(\thb)$ for most or all parameter points $\thb$, and 
locally D-optimal saturated designs at $\thb$ are known. Some models are presented 
in the following three examples. Moreover, the models in these examples satisfy condition (SI) or condition (${\rm GLM}^*$).
In a  fourth example  the model satisfies (${\rm GLM}^*$) and for almost all parameter points $\thb$   
the locally D-optimal saturated design at $\thb$ is unique and hence condition (SD)$(\thb)$ holds. Condition 
(SD$^*$)$(\thb)$ holds on a relevant subset of parameter points $\thb$ while on another subset  
(SD$^*$)$(\thb)$ does not hold, and for very special points $\thb$
(if included in $\Theta$) condition (SD)$(\thb)$ does not hold.

\vspace*{1ex}\noindent
{\em Example 1: Michaelis-Menten model.}\\
$p=2$, $\Theta\subseteq(\,0\,,\,\infty)^2$, ${\cal X}=[\,a\,,\,b\,]$ where $0\le a<b<\infty$, and
\begin{eqnarray*}
&&\mu(x,\theta)=\frac{\vartheta_1x}{\vartheta_2+x},\quad\mbox{where $\theta=(\vartheta_1,\vartheta_2)^\trp$;}\\
&&f_\theta(x)=\Bigl(\frac{\partial}{\partial\vartheta_1}\mu(x,\theta)\,,\,
\frac{\partial}{\partial\vartheta_2}\mu(x,\theta)\Bigr)^\trp=
\Bigl(\frac{x}{\vartheta_2+x}\,,\,-\frac{\vartheta_1x}{(\vartheta_2+x)^2}\Bigr)^\trp
\end{eqnarray*} 
for all $x\in[\,a\,,\,b\,]$ and $\theta=(\vartheta_1,\vartheta_2)^\trp\in\Theta$. 
For a given parameter point $\thb=(\overline{\vartheta}_1,\overline{\vartheta}_2)^\trp\in\Theta$, 
the unique locally D-optimal design at $\thb$
is the equally weighted two-point design \ $\xi^*_{\thb}$ supported by $x_1^*(\thb)$ and $b$, 
where $x_1^*(\thb)=\max\bigl\{\overline{\vartheta}_2b\big/(2\overline{\vartheta}_2+b)\,,\,a\bigr\}$,
see Bates and Watts \cite{Bates-Watts}, pp.~125-126. 
In fact, in that reference the design was shown to be the locally D-optimal saturated design at $\thb$.
Using the Kiefer-Wolfowitz equivalence theorem it can be checked that $\xi^*_{\thb}$ is locally D-optimal
at $\thb$. So the present model satisfies condition $({\rm SD}^*)(\thb)$ for all $\thb\in\Theta$. Moreover,
the model satisfies condition (SI), see \cite{FF-NG-RS-19}.

\vspace*{1ex}\noindent
{\em Example 2: Exponential decay model.}\\
$p=2$, $\Theta\subseteq(\,0\,,\,\infty)^2$, ${\cal X}=[\,a\,,\,b\,]$ where $0\le a<b<\infty$, and
\begin{eqnarray*}
&&\mu(x,\theta)=\vartheta_1\exp\bigl(-\vartheta_2x\bigr),\quad\mbox{where $\theta=(\vartheta_1,\vartheta_2)^\trp$;}\\
&&f_\theta(x)=\Bigl(\frac{\partial}{\partial\vartheta_1}\mu(x,\theta)\,,\,
\frac{\partial}{\partial\vartheta_2}\mu(x,\theta)\Bigr)^\trp
=\exp(-\vartheta_2x)\,\Bigl(\,1\,,\,-\vartheta_1x\,\Bigr)^\trp
\end{eqnarray*}
for all $x\in[\,a\,,\,b\,]$ and $\theta=(\vartheta_1,\vartheta_2)^\trp\in\Theta$.
For a given parameter point $\thb=(\overline{\vartheta}_1,\overline{\vartheta}_2)^\trp\in\Theta$, 
the unique locally D-optimal design at $\thb$
is the equally weighted two-point design \ $\xi^*_{\thb}$ supported by $a$ and $x_2^*(\thb)$, 
where $x_2^*(\thb)=\min\bigl\{a+1/\overline{\vartheta}_2\,,\,b\bigr\}$, see Box and Lucas \cite{Box-Lucas}, p.~85.
In fact, in that reference the design was shown to be the locally D-optimal saturated design at $\thb$.
Again, by the Kiefer-Wolfowitz equivalence theorem it can be verified that $\xi^*_{\thb}$ is locally D-optimal
at $\thb$. So the present model satisfies condition $({\rm SD}^*)(\thb)$ for all $\thb\in\Theta$. Moreover,
the model satisfies condition (SI), see \cite{FF-NG-RS-19}.

\vspace*{1ex}\noindent
{\em Example 3: Generalized linear models with binary response.}\\
Let $p=2$, $\Theta\subseteq\mathbb{R}^2$, ${\cal X}=[\,a\,,\,b\,]$ where $-\infty<a<b<\infty$.
Consider the class of generalized linear models given by
\[
\mu(x,\theta)=G\bigl(\vartheta_1+\vartheta_2 x)\quad \mbox{and}\quad
f_\theta(x)=\varphi\bigl(\vartheta_1+\vartheta_2 x\bigr)\,(1\,,\,x)^\trp,\quad \theta=(\vartheta_1,\vartheta_2)^\trp,
\]
where $G$ is a continuously differentiable 
distribution function on the real line with positive derivative $G'>0$, and 
\[
\varphi(u)=G'(u)\,\big/\,\sqrt{G(u)\bigl(1-G(u)\bigr)},\ \ u\in\mathbb{R}.
\]
The inverse function $G^{-1}$ is called the link function.
The models refer to binary response variables, and thus $\mu(x,\theta)$ equals the probability of a positive 
response at $x$. In particular, condition (${\rm GLM}^*$) is met, where $\psi(x,\theta)=\varphi(\vartheta_1+\vartheta_2 x)$.
Consider four particular members of that class of models:\\[.5ex]
(i) \ $G(u)\,=\,1\big/\bigl(1+\exp(-u)\bigr)$,\quad (logit link);\\[.5ex]
(ii) $G(u)\,=\,1-\exp\bigl\{-\exp(u)\bigr\}$,\quad (complementary log-log);\\[.5ex]
(iii) $G(u)\,=\,(2\pi)^{-1/2}\int_{-\infty}^u \exp(-t^2/2)\,{\rm d}t$,\quad (probit)\\[.5ex]
(iv) $G(u)\,=\, 1\big/\bigl(1+\exp(-u)\bigr)^m$, $m>0$ fixed,\quad (skewed logit).\\[1ex]
It was shown in Biedermann, Dette and Zhu \cite{Biedermann-Dette-Zhu} that, under each of models (i) to (iv), 
the locally D-optimal design at any given $\thb\in\Theta$ is unique and is an equally weighted two-point design.
Actually, in that paper a different parametrization of the models was employed and the results on local optimality
were obtained for a greater class of optimality criteria (Kiefer's criteria). 
For the D-criterion the locally D-optimal designs are equivariant under a parameter transformation, 
and therefore the results of \cite{Biedermann-Dette-Zhu} apply to the present models (i)--(iv), that is, 
the models satisfy condition $({\rm SD}^*)(\thb)$ for all $\thb\in\Theta$. 
For finding the support points of the locally D-optimal designs the results in Ford, Torsney and Wu 
\cite{Ford-Torsney-Wu}, Section 6, will be helpful. However, their derivations on p.~582 of the D-optimal 
saturated (two point) designs
are not conclusive. So we have included the result along with a proof in the appendix
as a supplement to this example (Appendix A.1).

\vspace*{1ex}\noindent
{\em Example 4: Poisson regression model with two covariates.}\\
Let $p=3$, $\Theta\subseteq\mathbb{R}\times(-\infty\,,\,0\,]^2$, ${\cal X}=[\,0\,,\,b_1]\times[\,0\,,\,b_2]$,
where $b_1>0$ and $b_2>0$. Consider a generalized linear model with Poisson distributed response variables,
\[
\mu(x,\theta)=\exp\bigl(\vartheta_0+\vartheta_1x_1+\vartheta_2x_2\bigr)\quad\mbox{and}\quad
f_\theta(x)= \exp\bigl\{{\textstyle\frac{1}{2}}(\vartheta_0+\vartheta_1x_1+\vartheta_2x_2)\bigr\}\,
(1\,,\,x_1\,,\,x_2)^\trp,
\]
where $x=(x_1,x_2)\in{\cal X}$ and $\theta=(\vartheta_0,\vartheta_1,\vartheta_2)^\trp\in\Theta$. 
In particular, condition (${\rm GLM}^*$) is met, where 
$\psi(x,\theta)= \exp\bigl\{{\textstyle\frac{1}{2}}(\vartheta_0+\vartheta_1x_1+\vartheta_2x_2)\bigr\}$,
$f(x)=(1\,,\,x_1\,,\,x_2)^\trp$, and $G(u)=\exp(u)$, $u\in\mathbb{R}$.  
Let $\thb=(\overline{\vartheta}_0,\overline{\vartheta}_1,\overline{\vartheta}_2)^\trp\in\Theta$ be given.
By our assumption on the parameter space  
the slope components $\overline{\vartheta}_1,\overline{\vartheta}_2$ are nonpositive. We consider three cases.\\[.5ex]
(i) $\overline{\vartheta}_1<0$, $\overline{\vartheta}_2<0$; \ (ii) $\overline{\vartheta}_1<0$, $\overline{\vartheta}_2=0$;
 \ (iii) $\overline{\vartheta}_1=\overline{\vartheta}_2=0$.\\[.5ex]
By standard arguments, the problem of finding a locally D-optimal saturated design at $\thb$ can equivalently be transformed
to that of finding a D-optimal saturated design for the linear regression model given by
$f_0(z)$, $z=(z_1,z_2)\in {\cal Z}=[0,c_1]\times[0,c_2]$, where\\
in case (i): $z_j=|\overline{\vartheta}_j|x_j$, $c_j=|\overline{\vartheta}_j|b_j$,  $j=1,2$, and 
$f_0(z)\,=\,\exp\bigl\{{\textstyle-\frac{1}{2}}(z_1+z_2)\bigr\}\,(1\,,\,z_1\,,\,z_2)^\trp$;\\
in case (ii): $z_1=|\overline{\vartheta}_1|x_1$, $z_2=x_2$, $c_1=|\overline{\vartheta}_1|b_1$, $c_2=b_2$, and 
$f_0(z)\,=\,\exp\bigl\{{\textstyle-\frac{1}{2}}z_1\bigr\}\,(1\,,\,z_1\,,\,z_2)^\trp$;\\
in case (iii): $z_j=x_j$, $c_j=b_j$, $j=1,2$, and $f_0(z)=(1\,,\,z_1\,,\,z_2)^\trp$.\\  
Lemma \ref{lemA-1} in Appendix A.2 yields the D-optimal saturated designs in terms of the $z$-variable,
which are easily transformed back to the locally D-optimal saturated designs in the original model. 
In case (i) the locally D-optimal saturated design is unique and hence condition (SD)$(\thb)$ holds;
in cases (ii) and (iii) there are infinitly many  locally D-optimal saturated designs and, as it is easily seen,
their information matrices vary, hence condition (SD)$(\thb)$ does not hold. Furthermore, in case (i) the following 
holds (see Lemma \ref{lemA-3} in Appendix A.2). 
If $|\overline{\vartheta}_j|\ge2/b_j$ for $j=1,2$ then the locally D-optimal saturated design is locally D-optimal
and hence condition $({\rm SD}^*)(\thb)$ holds. On the other hand, if  
$|\overline{\vartheta}_1|$ and $|\overline{\vartheta}_2|$ are small in the sense that  $|\overline{\vartheta}_j|\le2/b_j$ for $j=1,2$ 
and $\bigl(1+\exp(-|\overline{\vartheta}_1|\,b_1\bigr)\,\bigl(1+\exp(-|\overline{\vartheta}_2|\,b_2\bigr)\,>2$, then
the locally D-optimal saturated design      
is not locally D-optimal and hence condition $({\rm SD}^*)(\thb)$ does not hold.

Poisson models with two or more covariates
were considered by Russell et al.~\cite{Russell-et-al} and more general results on locally D-optimal designs were
obtained. In their Remark 3 on p.~724 a result on locally D-optimal saturated designs covering case (i) of the present model 
was stated but no proof was given. We give a proof in Appendix A.2 for the present situation of two covariates.
\eop

\section{Adaptive $p$-step-ahead algorithm}
\setcounter{equation}{0}
Let $\mathbb{N}$ denote the set of all positive integers.
By $\delta[x]$, for any $x\in{\cal X}$, 
we denote the one-point distribution on ${\cal X}$ concentrated at the point $x$.
The adaptive algorithm described next generates iteratively (in batches of size $p$) a sequence of 
design points. For each batch of design points the responses are observed, and
the parameter estimate is updated based on all design points and responses obtained so far. 
The estimate is used for choosing the next batch of design points, and so on. 
Along with the sequences of design points and response values, 
a sequence of designs and a sequence of parameter estimates emerge.      

\vspace*{2ex}\noindent
\underline{Algorithm}
\begin{itemize}
\item[(o)] {\em Initialization ($k=1$):} A number $n_1\in\mathbb{N}$ and design points $x_1,\ldots,x_{n_1}\in{\cal X}$ 
are chosen forming the initial design $\xi_1=\frac{1}{n_1}\sum_{i=1}^{n_1}\delta[x_i]$.
Observations $y_1,\ldots,y_{n_1}$ of responses at the design points $x_1,\ldots,x_{n_1}$, respectively, are taken.
Based on the current data a parameter estimate $\theta_1\in\Theta$ is computed,
\[
\theta_1=\widehat{\theta}_1(x_1,y_1,\ldots,x_{n_1},y_{n_1}).
\]
\item[(i)] {\em Iteration:} Let $k\ge1$ and $n_k=n_1+(k-1)p$, let the current data be given by the points $x_1,\ldots,x_{n_k}\in{\cal X}$
forming the current design $\xi_k=\frac{1}{n_k}\sum_{i=1}^{n_k}\delta[x_i]$, and by the observed responses $y_1,\ldots,y_{n_k}$
at $x_1,\ldots,x_{n_k}$, respectively, and let $\theta_k=\widehat{\theta}_k(x_1,y_1,\ldots,x_{n_k},y_{n_k})$ be the current parameter estimate  
on the basis of the current data. Then, a batch of $p$ design points $x_{n_k+1},\ldots,x_{n_k+p}\in{\cal X}$ is chosen such that
\begin{equation}
\bigl(\det\bigl[f_{\theta_k}(x_{n_k+1}),\ldots,f_{\theta_k}(x_{n_k+p})\bigr]\bigr)^2\,=\,
\max_{z_1,\ldots,z_p\in{\cal X}}\bigl(\det\bigl[f_{\theta_k}(z_1),\ldots,f_{\theta_k}(z_p)\bigr]\bigr)^2.
\label{eq2-1}
\end{equation}

Observations $y_{n_k+1},\ldots,y_{n_k+p}$ of responses at $x_{n_k+1},\ldots,x_{n_k+p}$, respectively, are taken and, 
based on the augmented data, a new parameter estimate $\theta_{k+1}\in\Theta$  is computed,
\[
\theta_{k+1}\,=\,\widehat{\theta}_{k+1}(x_1,y_1,\ldots,x_{n_k+p},y_{n_k+p}).
\]
Set $n_{k+1}=n_k+p$ and $\xi_{k+1}=\bigl(1/n_{k+1}\bigr)\sum_{i=1}^{n_{k+1}}\delta[x_i]$.
Iteration step (i) is repeated with $k$ replaced by $k+1$. 
\eop
\end{itemize}

\noindent
{\em Remarks.}\\[.5ex] 
1. Obviously, in the iteration step (i) we have $\xi_{k+1}=(n_k/n_{k+1})\,\xi_k\,+\,(p/n_{k+1})\,\eta_k$, where\\
$\eta_k\,=\,(1/p)\sum_{j=1}^p\delta[x_{n_k+j}]$, and by (\ref{eq2-1}) $\eta_k$
is a locally D-optimal saturated design at $\theta_k$.\\
2. For the initial design of the algorithm, $\xi_1=(1/n_1)\sum_{i=1}^{n_1}\delta[x_i]$, 
the number $n_1$ of points (and the points themselves) may be arbitrary. In practice, one 
might prefer some saturated design and thus $n_1=p$. The choice $n_1=p$ will also simplify some 
theoretical derivations in Sections 5 and 6. In fact, in our proofs of the theorems 
we will assume $n_1=p$ to cut down the technical effort. However, the results hold for any choice of $n_1$.\\[.5ex]
3.  The adaptive Wynn algorithm studied in \cite{FF-NG-RS-18} and \cite{FF-NG-RS-19}
requires that the initial design $\xi_1$ is such that its information matrix $M(\xi_1,\theta)$
is non-singular for all $\theta\in\Theta$, which implies that all subsequently generated designs $\xi_k$, $k\ge2$,
have that property as well. The iteration rule of the adaptive Wynn algorithm is given by  
\[
x_{k+1}\,=\,\arg\max_{x\in{\cal X}}f_{\theta_k}^\trp(x)\, M^{-1}(\xi_k,\theta_k)\,f_{\theta_k}(x).
\]
In the (nearly) trivial case $p=1$ this becomes $x_{k+1}=\arg\max_{x\in{\cal X}}\bigl(f_{\theta_k}(x)\bigr)^2$ 
which coincides with the iteration rule of the present $p$-step-ahead algorithm 
in case $p=1$. So, for $p=1$, the present algorithm coincides with the adaptive Wynn algorithm.
Note also, that for $p=1$ condition (SD$^*$)$(\thb)$ holds for any $\thb\in\Theta$, since
a locally D-optimal design at $\thb$ is given by the one-point design $\delta[x^*_{\thb}]$,
where $x^*_{\thb}=\arg\max_{x\in{\cal X}}\bigl(f_{\thb}(x)\bigr)^2$. 
\eop     

The algorithm uses observations of responses which are values of 
random variables (response variables). So the generated sequences $x_i$, $y_i$ ($i\in\mathbb{N}$) and
$\xi_k$, $\theta_k$ ($k\in\mathbb{N}$) are random and should be viewed as paths of corresponding sequences of random variables.
This will be modeled appropriately in Sections 5 and 6. 
In Section 4 we focus on some properties of the algorithm which do not require a specific stochastic model.
The proofs of the results have been transferred to the appendix (parts A.3 and A.4).

\section{Some basic properties of the algorithm}
\setcounter{equation}{0}
The Euclidean norm in $\mathbb{R}^p$ is given by $\Vert a\Vert=(a^\trp a)^{1/2}$. The Frobenius norm in the space 
$\mathbb{R}^{p\times p}$ of all 
$p\times p$ matrices is given by $\Vert A\Vert = \Bigl(\sum_{i,j=1}^pa_{ij}^2\Bigr)^{1/2}$ for $A=(a_{ij})_{1\le i,j\le p}$.
For a symmetric $p\times p$ matrix $A$ the smallest eigenvalue of $A$ is denoted by $\lambda_{\rm\scriptsize min}(A)$.
The distance function in the (compact) metric space $\Theta$ is denoted by ${\rm d}_\Theta$, and the
set of all designs on ${\cal X}$ is denoted by $\Xi$.
We start with an auxiliary lemma which does not specifically refer to the algorithm. 

\begin{lemma}\label{lem3-1} 
Let $\rho_k, \tau_k\in\Theta$, $k\in\mathbb{N}$, 
be two sequences of parameter points such that\\ 
$\lim_{k\to\infty}{\rm d}_{\Theta}(\rho_k,\tau_k)\,=0$.
Then
\[
\lim_{k\to\infty}\Bigl(\sup_{\xi\in\Xi}\big\Vert M(\xi,\rho_k)\,-\,
M(\xi,\tau_k)\big\Vert\Bigr)\,=\,0.
\]
As a consequence, if $\Phi$ is a real-valued continuous function on the set of all nonnegative definite 
$p\times p$ matrices, then
\[
\lim_{k\to\infty}\Bigl(
\sup_{\xi\in\Xi}\big|\Phi\bigl(M(\xi,\rho_k)\bigr)-\Phi\bigl(M(\xi,\tau_k)\bigr)
\big|\Bigr)\,=\,0.
\]
\end{lemma}

For $B\in\mathbb{R}^{p\times p}$ and $\emptyset\not={\cal A}\subseteq\mathbb{R}^{p\times p}$, we denote
by ${\rm dist}(B,{\cal A})$ the distance of the point $B$ and the set ${\cal A}$, that is, 
${\rm dist}(B,{\cal A})=\inf_{A\in {\cal A}}\Vert B-A\Vert$.
As it is well-known, the function $B\mapsto{\rm dist}(B,{\cal A})$ on $\mathbb{R}^{p\times p}$ is continuous, and
if the set ${\cal A}$ is convex then this function is convex. For any nonempty subset ${\cal A}\subseteq\mathbb{R}^{p\times p}$
we denote by ${\rm Conv}\,{\cal A}$ the convex hull of ${\cal A}$, that is,
\[
{\rm Conv}\,{\cal A}\,=\,\Bigl\{\sum_{i=1}^r\alpha_iA_i\,:\,\alpha_i\ge0,\  A_i\in{\cal A}\ (1\le i\le r),
\ \sum_{i=1}^r\alpha_i=1, \ r\in\mathbb{N}\,\Bigr\}.
\]
As a particular set ${\cal A}$ we consider the set of information matrices at $\thb$ of all 
locally D-optimal saturated designs at $\thb$, for a given parameter point $\thb\in\Theta$. We denote
\[
{\cal M}_{{\rm\scriptsize s}*}(\thb)\,=\,
\bigl\{M(\eta^*,\thb)\,:\, \mbox{$\eta^*$ is a locally D-optimal saturated design at $\thb$}\bigr\}.
\]
In the following lemma an arbitrary path of the $p$-step-ahead algorithm is considered 
yielding a sequence $\xi_k$ of designs and a sequence $\theta_k$ of
parameter estimates.     

\begin{lemma}\label{lem3-2}
If $\lim_{k\to\infty}\theta_k=\overline{\theta}$ 
for some $\overline{\theta}\in\Theta$,
then for every sequence $\theta_k'\in\Theta$, $k\in\mathbb{N}$, such that 
$\lim_{k\to\infty}\theta_k'=\overline{\theta}$ one has
\[
{\rm dist}\Bigl(M(\xi_k,\theta_k')\,,\,{\rm Conv}\,{\cal M}_{{\rm\scriptsize s}*}(\overline{\theta})\Bigr)\,\longrightarrow\,0
\ \mbox{ as $k\to\infty$.}
\]
Under condition $({\rm SD})(\overline{\theta})$ the latter convergence is the same as 
\[  
\lim_{k\to\infty}M(\xi_k,\theta_k')\,=\,M_{{\rm\scriptsize s}*}(\overline{\theta}),
\]
with $M_{{\rm\scriptsize s}*}(\overline{\theta})$ according to condition $({\rm SD})(\overline{\theta})$.
\end{lemma}

We denote the distance function in the (compact) metric space ${\cal X}$ by ${\rm d}_{\cal X}$. Again, 
we consider any path of the $p$-step-ahead algorithm, and now we focus on the sequences  
$x_i$ ($i\in\mathbb{N}$) and $\xi_k$ ($k\in\mathbb{N}$) of design points and designs, respectively.

\begin{lemma}\label{lem3-3}\quad\\
(i) \ There exists a constant $\Delta_0>0$ such that
\[
{\rm d}_{{\cal X}}(x_{n_k+\ell},x_{n_k+m})\,\ge\Delta_0\ \mbox{ for all $1\le\ell<m\le p$ and all $k\ge1$.}
\] 
(ii) \ Under condition {\rm(GLM)}, there exists a constant  $\varepsilon_0>0$ such that
\[
\xi_k\Bigl(\bigl\{x\in{\cal X}\,:\,|a^\trp f(x)|\ge\varepsilon_0\bigr\}\Bigr)\,\ge\,(k-1)/n_k \ 
\mbox{ for all $a\in\mathbb{R}^p$, $\Vert a\Vert=1$, and all $k\ge1$.}
\]
\end{lemma}

\noindent
{\em Remark.} 
The constants $\Delta_0$ and $\varepsilon_0$ constructed in the proof of Lemma \ref{lem3-3} (see Appendix A.3)
depend only on the family $f_\theta$, $\theta\in\Theta$, but they do not depend on the particular path generated by the $p$-step-ahead algorithm.
So $\Delta_0$ and $\varepsilon_0$ in the lemma can be chosen simultaneously for all possible paths of the algorithm.
\eop

A desirable property of a sequence of estimators of $\theta$ is strong consistency, that is, almost sure convergence to the true 
parameter point $\overline{\theta}$. 
For a sequence of random variables $W_k$, $k\in\mathbb{N}$, and a random variable $W$ defined on some probability space
with values in some metric space, the notation $W_k\asto W$ stands for almost sure convergence of $W_k$ to $W$ 
as $k\to\infty$. Under the assumption that the estimators $\widehat{\theta}_k$, $k\in\mathbb{N}$, employed by the algorithm are
strongly consistent, aymptotic properties of the designs $\xi_k$, $k\in\mathbb{N}$, generated by the algorithm are stated as a corollary
below. A desirable property is ``asymptotic local D-optimality at $\thb$ (almost surely)'', that is, $\det M(\xi_k,\thb)\,\asto d_*(\thb)$  where
$d_*(\thb)$ denotes the maximum value of $\det M(\xi,\thb)$ over all designs $\xi$.
It is not difficult to show that asymptotic local D-optimality at $\thb$ of the sequence $\xi_k$ is equivalent to \    
$M(\xi_k,\overline{\theta}) \asto M_*(\overline{\theta})$, where $M_*(\overline{\theta})$ is the unique information matrix at $\overline{\theta}$ of a   
locally D-optimal design at $\overline{\theta}$. Since the concept of the $p$-step-ahead algorithm is based on locally D-optimal saturated 
designs, one cannot expect asymptotic local D-optimality at $\thb$ (a{.}s{.}) of the design sequence $\xi_k$ in general, unless condition     
$({\rm SD}^*)(\overline{\theta})$ holds. The following corollary is a fairly direct consequence of Lemmas \ref{lem3-1} and \ref{lem3-2},
and we thus state it without a proof. Recall notations  ${\cal M}_{{\rm\scriptsize s}*}(\overline{\theta})$ for the set of 
information matrices at $\thb$ of all locally D-optimal saturated designs at $\thb$ and, in the case that condition $({\rm SD})(\overline{\theta})$
holds, $\Ms(\overline{\theta})$ for the unique element of ${\cal M}_{{\rm\scriptsize s}*}(\overline{\theta})$. Furthermore, let
$d_{{\rm\scriptsize s}*}(\thb)$ be the maximum value of $\det M(\eta,\thb)$ over all saturated designs $\eta$.

\begin{corollary}\label{cor3-1}
Assume that the sequence of adaptive estimators $\widehat{\theta}_k$,
$k\in\mathbb{N}$, employed by the $p$-step-ahead algorithm is strongly consistent, that is, $\widehat{\theta}_k\asto\thb$ 
where $\thb\in\Theta$ is the true parameter point. 
Then, for the sequence of designs $\xi_k$, $k\in\mathbb{N}$, generated by the algorithm one has:\\[.3ex]
{\rm(i)} \ ${\rm dist}\Bigl(M(\xi_k,\thb)\,,\,{\rm Conv}\,{\cal M}_{{\rm\scriptsize s}*}(\overline{\theta})\Bigr)\,\asto 0$ \ and hence
$\liminf_{k\to\infty} \det M(\xi_k,\thb)\,\ge d_{{\rm\scriptsize s}*}(\thb)$ a.s.\\[.3ex]
{\rm(ii)} If condition $({\rm SD})(\overline{\theta})$ holds, then \ $M(\xi_k,\thb) \asto \Ms(\overline{\theta})$ \ and \ 
$\det M(\xi_k,\thb) \asto d_{{\rm\scriptsize s}*}(\thb)$.\\[.6ex]
{\rm(iii)} If condition (SD$^*$)$(\overline{\theta})$ holds, then the designs $\xi_k$ are asymptotically locally D-optimal at $\thb$ (a{.}s{.}), that is,
\ $M(\xi_k,\overline{\theta}) \asto M_*(\overline{\theta})$ \ and \ $\det M(\xi_k,\thb)\,\asto d_*(\thb)$.  
\end{corollary}

In Sections 5 and 6 we will show that adaptive least squares estimators and maximum likelihood estimators in the $p$-step-ahead algorithm are 
strongly consistent, under appropriate assumptions. 
In particular, the models in Examples 1 to 4 of Section 2 will be covered with adaptive least squares estimation 
in the Michaelis-Menten model (Example 1) and the exponential decay model (Example 2), 
and with adaptive maximum likelihood estimation in the generalized linear models of Examples 3 and 4. So for those models, 
when the algorithm employs least squares estimators and maximum likelihood estimators, respectively, by Corollary \ref{cor3-1}
the adaptive design sequence $\xi_k$ generated by the algorithm is asymptotically locally D-optimal at $\thb$ (a{.}s{.}) 
for any true parameter point $\thb\in\Theta$ in Examples 1 to 3, 
and for any true parameter point $\thb\in\Theta'\subseteq\Theta$ in Example 4 with a relevant subset of $\Theta'$.

\section{Adaptive least squares estimators}
\setcounter{equation}{0}
In this section and in the next, we will examine the asymptotic properties (strong consistency and asymptotic normality) of
adaptive least squares and adaptive maximum likelihood estimators in the  $p$-step-ahead algorithm. 
To this end, appropriate stochastic models for the algorithm will be employed. 

Let $X_i$ and $Y_i$, $i\in\mathbb{N}$, be two sequences of random variables defined on a common probability space 
$(\Omega,{\cal F},\mathbb{P}_{\overline{\theta}})$ where $\overline{\theta}\in\Theta$ denotes the true parameter point (which is unknown).
The random variables $X_i$ have their values in  ${\cal X}$ and the $Y_i$ are real valued. 
A run of the algorithm generates paths $x_i$ and $y_i$, $i\in\mathbb{N}$,
of the sequences $X_i$ and $Y_i$, respectively. An appropriate adaptive version of a regression model is stated 
by the following two assumptions (a1) and (a2), cf.~\cite{FF-NG-RS-19}, Section 3. 
Later, some further strengthening assumptions will be added. 

\begin{itemize}
\item[(a1)] 
Let a nondecreasing sequence \ 
${\cal F}_0\subseteq{\cal F}_1\subseteq\,\ldots\,\subseteq{\cal F}_k\subseteq\,\ldots$ \ 
of sub-sigma-fields of ${\cal F}$ be given such that  for each $k\in\mathbb{N}$ the multivariate random variable\\ 
$\BM{X}_k=(X_{n_{k-1}+1},\ldots,X_{n_k})$ \ 
is ${\cal F}_{k-1}$-measurable, 
and the multivariate random variable\\ 
$\BM{Y}_k=(Y_{n_{k-1}+1},\ldots,Y_{n_k})^\trp$ \ is ${\cal F}_k$-measurable.
Here we define $n_0:=0$.
\item[(a2)]  
$Y_i\,=\,\mu(X_i,\overline{\theta})\,+\,e_i$ \ for all $i\in\mathbb{N}$ with real-valued 
square integrable random errors $e_i$, $i\in\mathbb{N}$, such that the multivariate error variables 
$\BM{e}_k:=(e_{n_{k-1}+1},\ldots,e_{n_k})^\trp$, $k\in\mathbb{N}$, satisfy: 
${\rm E}\bigl(\BM{e}_k\,\big|\,{\cal F}_{k-1}\bigr)\,=0 \ \mbox{\,a.s.}$\  for all $k\in\mathbb{N}$, and \ 
$\sup_{k\in\mathbb{N}}{\rm E}\bigl(\Vert\BM{e}_k\Vert^2\,\big|\,{\cal F}_{k-1}\bigr)\,<\infty \ \ \mbox{a.s.}$
\end{itemize}

Since $n_k=n_1+(k-1)p$ for all $k\ge1$, the dimensions of the multivariate random variables $\BM{X}_k$, $\BM{Y}_k$, and
$\BM{e}_k$ introduced in (a1) and (a2)
are given by $n_k-n_{k-1}=p$ for all $k\ge2$ and $n_1-n_0=n_1$. In the proofs of consistency and asymptotic normality we will
restrict to the case $n_1=p$. 
  
The adaptive least squares estimators (adaptive LSEs)
$\widehat{\theta}_k^{({\scriptsize\rm LS})}=\widehat{\theta}_k^{({\scriptsize\rm LS})}(X_1,Y_1,\ldots,X_{n_k},Y_{n_k})$, $k\ge1$,
are defined pathwise by
\[
\widehat{\theta}_k^{({\rm\scriptsize LS})}(x_1,y_1,\ldots,x_{n_k},y_{n_k})\,=\,\arg\min_{\theta\in\Theta}  
\sum_{i=1}^{n_k}\bigl(y_i-\mu(x_i,\theta)\bigr)^2.
\]
Note that we do not generally assume that 
the adaptive estimators employed by the algorithm, 
$\widehat{\theta}_k=\widehat{\theta}_k(X_1,Y_1,\ldots,X_{n_k},Y_{n_k})$, $k\ge1$, 
are given by the adaptive LSEs. 

Under condition (SI) of `saturated identifiability' or, alternatively, condition $({\rm GLM}^*)$ of `generalized linear model',
strong consistency of the adaptive LSEs is shown by the next result. 
Note that the adaptive estimators $\widehat{\theta}_k$ employed by the algorithm
may be arbitrary.

\begin{theorem}\label{theo4-1}
Assume model {\rm(a1)}, {\rm(a2)}, and assume one of conditions {\rm(SI)} or $({\rm GLM}^*)$. Then: 
\ $\widehat{\theta}_k^{({\rm\scriptsize LS})}\asto\overline{\theta}$.
\end{theorem}

For achieving asymptotic normality further conditions are needed. Firstly, the basic conditions (assumed throughout)
(b1)-(b4) are augmented by the `gradient condition' (b5) on the family of functions $f_\theta$, $\theta\in\Theta$,  and the mean response $\mu$. 
\begin{itemize}
\item[(b5)]   
$\Theta\subseteq\mathbb{R}^p$ (endowed with the usual Euclidean metric), 
${\rm int}(\Theta)\not=\emptyset$, where ${\rm int}(\Theta)$ denotes the interior of $\Theta$ as a subset of $\mathbb{R}^p$,
the function \ $\theta\mapsto\mu(x,\theta)$ \ is twice differentiable on the interior of $\Theta$
for each fixed $x\in{\cal X}$, with gradients and Hessian matrices denoted by 
$\nabla \mu(x,\theta)=
\Bigl(\frac{\partial}{\partial\vartheta_1}\mu(x,\theta),\ldots,\frac{\partial}{\partial\vartheta_p}\mu(x,\theta)\Bigr)^\trp$
 \ and \ 
$\nabla^2\mu(x,\theta)=\Bigl(\frac{\partial^2}{\partial\vartheta_i\partial\vartheta_j}\mu(x,\theta)\Bigr)_{1\le i,j\le p}$, respectively, where 
$\theta=(\vartheta_1,\ldots,\vartheta_p)^\trp$. 
The functions $(x,\theta)\mapsto\nabla\mu(x,\theta)$ and $(x,\theta)\mapsto\nabla^2\mu(x,\theta)$ 
are continuous on ${\cal X}\times {\rm int}(\Theta)$, and
\[
f_\theta(x)\,=\,\nabla \mu(x,\theta)\quad\mbox{for all $x\in{\cal X}$ and all $\theta\in{\rm int}(\Theta)$.}
\] 
\end{itemize}

Two additional conditions (L) and (AH) on the error variables 
of model (a1), (a2) are imposed, where `L' stands for `Lindeberg' and `AH' for `asymptotic homogeneity'. 
For an event $A$ in the underlying probability space we denote by $\Ifkt(A)$ the dichotomous random variable which 
yields the value $1$ if the event $A$ occurs, and yields the value $0$ otherwise.

\begin{itemize}
\item[(L)]
$\displaystyle\frac{1}{k}\sum_{j=1}^k 
{\rm E}\Bigl(\Vert \BM{e}_j\Vert^2\Ifkt\bigl(\Vert \BM{e}_j\Vert>\varepsilon\sqrt{k}\bigr)\,\Big|{\cal F}_{j-1}\Bigr)
\,\asto0$ \ \ for all $\varepsilon>0$.
\item[(AH)]
$\displaystyle{\rm E}\bigl(\BM{e}_k\BM{e}_k^\trp\big|{\cal F}_{k-1}\bigr)\asto\sigma^2(\overline{\theta})\,I_p$ \ 
 \  for some positive real constant $\sigma^2(\overline{\theta})$, where
$I_p$ denotes the $(p\times p)$ identity matrix.
\end{itemize}

Each of the following two conditions (L') and (L'') implies (L), 
which can be seen by similar arguments as in \cite{FF-NG-RS-19}, Section 3.

\begin{itemize}
\item[(L')] 
$\sup_{k\in\mathbb{N}}{\rm E}\bigl(\Vert\BM{e}_k\Vert^\alpha\big|{\cal F}_{k-1}\bigr)\,<\,\infty$ \ a.s. \ for some $\alpha>2$.
\item[(L'')] 
The random variables $\BM{e}_k$, $k\ge2$, are identically distributed, and \ 
$\BM{e}_k$, ${\cal F}_{k-1}$ \ are independent for each $k\ge2$.
\end{itemize}

The $m$-dimensional normal distribution with expectation $0$ and covariance matrix $C$ is denoted by 
${\rm N}(0,C)$, where $C$ is a positive definite $m\times m$ matrix. 
For a sequence $W_k$ of $\mathbb{R}^m$-valued random variables, 
convergence in distribution of $W_k$ (as $k\to \infty$) to the $m$-dimensional normal distribution ${\rm N}(0,C)$ 
is abbreviated by $W_k\dto{\rm N}(0,C)$.   
  
\begin{theorem}\quad\label{theo4-2}\\
Assume model {\rm(a1)}, {\rm(a2)}, and assume conditions  {\rm(b5)}, {\rm(L)}, {\rm(AH)}, 
and $({\rm SD})(\overline{\theta})$. Moreover, assume
that the sequence $\widehat{\theta}_k$ of adaptive estimators employed by the algorithm and the sequence of adaptive LSEs
$\widehat{\theta}_k^{\rm\scriptsize (LS)}$ are both strongly consistent, that is,
$\widehat{\theta}_k\asto\overline{\theta}$ and $\widehat{\theta}_k^{\rm\scriptsize (LS)}\asto\thb$,
and let $\overline{\theta}\in{\rm int}(\Theta)$. Then: 
\[
\sqrt{n_k}\,\bigl(\widehat{\theta}_k^{\rm\scriptsize (LS)}-\overline{\theta}\bigr)
\,\dto\,{\rm N}\bigl(0,\sigma^2(\thb)\,\Ms^{-1}(\thb)\bigr),
\]
with $\Ms(\thb)$ according to condition $({\rm SD})(\overline{\theta})$.
\end{theorem}

\section{Adaptive maximum likelihood estimators}
In this section we consider an adaptive version of a generalized linear model. 
Let a one-parameter exponential family $P_\tau$, $\tau\in J$ be given, where $J\subseteq\mathbb{R}$ is an open interval
and $\tau$  is the canonical parameter. The $P_\tau$ are probability distributions on the Borel sigma-field of the real line with
densities w.r.t. some Borel-measure $\nu$,
\begin{equation}
p_\tau(y)\,=\,K(y)\,\exp\Bigl(\tau\,y\,-\,b(\tau)\Bigr),\ \ y\in\mathbb{R},\ \tau\in J,
\label{eq5-1}
\end{equation}
where $K$ is a nonnegative measurable function on $\mathbb{R}$ and $b$ is a real-valued function on $J$.
The function $b$ is infinitely differentiable, and for
its first and second derivatives one has $b'(\tau)={\rm E}_\tau(Y)$ and $b''(\tau)={\rm Var}_\tau(Y)>0$, the expectation and the variance of $P_\tau$,
respectively, see Fahrmeir and Kaufmann \cite{Fahrmeir-Kaufmann}, Section 2. 
In particular, the first derivative $b'$ is a smooth and strictly increasing function and hence a bijection,
$b'\,:\,J\longrightarrow b'(J)$, where the image $b'(J)$ is an open interval of the real line
and equals the set of expectations $\bigl\{E_\tau(Y)\,:\,\tau\in J\bigr\}$. 
Condition (${\rm GLM}^*$) is assumed where the scalar-valued function $\psi(x,\theta)$ in (GLM)
is given by
\begin{eqnarray}
&&\psi(x,\theta)=\varphi\bigl(f^\trp(x)\,\theta\bigr) \ \ \mbox{ for all $(x,\theta)\in{\cal X}\times\Theta$,}\label{eq-a3'}\\
&&\mbox{where }\ \varphi(u)\,=\,G'(u)\Big/\sqrt{b''\Bigl((b')^{-1}\bigl(G(u)\bigr)\Bigr)}\ \mbox{ for all $u\in I$,}\nonumber
\end{eqnarray}
and where it is assumed that $G(I)\subseteq b'(J)$.

As in Section 5 let $X_i$ and $Y_i$, $i\in\mathbb{N}$, be two sequences of random variables defined on a probability space
$(\Omega,{\cal F},\mathbb{P}_{\overline{\theta}})$ and with values in ${\cal X}$ and $\mathbb{R}$, respectively,
where $\overline{\theta}\in\Theta$ denotes the true (but unknown) parameter point. 
The stochastic model for the adaptive algorithm is given by assumption (a1) 
from Section 5 plus the following (a2'), which is stronger than (a2) from Section 5.
Recall the multivariate random variables $\BM{X}_k=(X_{n_{k-1}+1},\ldots,X_{n_k})$ 
and  $\BM{Y}_k=(Y_{n_{k-1}+1},\ldots,Y_{n_k})^\trp$, $k\in\mathbb{N}$.

\begin{itemize}
\item[(a2')]
For each $k\in\mathbb{N}$ the conditional distribution of $\BM{Y}_k$ given ${\cal F}_{k-1}$
is equal to the product of the distributions $P_{\overline{\tau}_i}$, $n_{k-1}+1\le i\le n_k$, where 
$\overline{\tau}_i=(b')^{-1}\Bigl(G\bigl(f^\trp(X_i)\,\overline{\theta}\bigr)\Bigr)$.
\end{itemize}

To interprete the random variables $\overline{\tau}_i$ in (a2') we note that for any $x\in{\cal X}$ the
parameter value $\overline{\tau}(x)= (b')^{-1}\Bigl(G\bigl(f^\trp(x)\,\overline{\theta}\bigr)\Bigr)$ selects that distribution
$P_{\overline{\tau}(x)}$ from the exponential family whose expectation equals $G\bigl(f^\trp(x)\,\overline{\theta}\bigr)$,
according to condition (${\rm GLM}^*$). 
Note that for the canonical link, that is $I=J$ and $G=b'$,  formulas simplify to 
$\tau(x)=f^\trp(x)\,\thb$, $\overline{\tau}_i=f^\trp(X_i)\,\thb$, and $\varphi(u)=\sqrt{b''(u)}$.
Note further that (${\rm GLM}^*$) together with (\ref{eq-a3'}) 
ensures that the information matrices from (\ref{eq1-1})
yield the Fisher information matrices, see Atkinson and Woods \cite{Atkinson-Woods}, p.~473, see also 
Fahrmeir and Kaufmann \cite{Fahrmeir-Kaufmann}, p.~347.

\vspace*{1ex}\noindent
{\em Example 3 (continued).}\\
Consider the class of generalized linear models with binary response from Example 3 in Section 2.
The family of binomial-$(1,\pi)$-distributions (where $0<\pi<1$) rewrites in canonical form (\ref{eq5-1})
with canonical parameter $\tau=\log\bigl(\pi/(1-\pi)\bigr)\in\mathbb{R}$ and 
$b(\tau)=\log\bigl(1+\exp(\tau)\bigr)$.
The densities refer to the two-point Borel measure $\nu=\delta[0]+\delta[1]$, and $K(y)=1$ if $y\in\{0,1\}$,
and $K(y)=0$ else. By straightforward calculation,
\[
b'(\tau)=\exp(\tau)\big/\bigl(1+\exp(\tau)\bigr),\ \  (b')^{-1}(\pi)=\log\bigl(\pi/(1-\pi)\bigr),\ \mbox{ and }\ 
b''(\tau)= \exp(\tau)\big/\bigl(1+\exp(\tau)\bigr)^2.
\]
Hence 
\[
b''\Bigl((b')^{-1}\bigl(G(u)\bigr)\Bigr)\,=\,G(u)\,\bigl(1-G(u)\bigr),
\]
which shows that the function $\varphi$ employed in Example 3 of Section 2 corresponds to (\ref{eq-a3'}).
Note that the logit model (i) of the example employs the canonical link, 
$G=b'$, and hence for this model 
$\varphi(u)=\sqrt{b''(u)}=\exp(u/2)\big/\bigl(1+\exp(u)\bigr)$. 
\eop

As in \cite{FF-NG-RS-18}, Section 3, one concludes from (a1), (a2') that the joint log-likelihood of  $X_1,Y_1,\ldots,X_{n_k},Y_{n_k}$
(up to an additive term not depending on $\theta$) is given by 
\begin{eqnarray}
&&L_{n_k}(\theta)\,=\,\sum_{i=1}^{n_k}\Bigl(\log\bigl(K(Y_i)\bigr) + \tau_i(\theta)\,Y_i - b\bigl(\tau_i(\theta)\bigr)\Bigr), \label{eq5-2}\\
&&\mbox{ where }\ \tau_i(\theta)\,=\,(b')^{-1}\Bigl(G\bigl(f^\trp(X_i)\,\theta\bigr)\Bigr).\label{eq5-2a}
\end{eqnarray} 
The adaptive maximum likelihood estimator $\MLk=\MLk(X_1,Y_1,\ldots,X_{n_k},Y_{n_k})$ maximizes $L_{n_k}(\theta)$ over $\theta\in\Theta$. 
Its strong consistency is shown by the next result. Note that the adaptive estimators $\widehat{\theta}_k$
employed by the algorithm may be arbitrary.  

\begin{theorem}\label{theo5-1}
Assume {\rm(a1)}, {\rm(a2')}, and $({\rm GLM}^*)$ with (\ref{eq-a3'}). Then \ $\widehat{\theta}_k^{(\rm\scriptsize ML)}\asto \overline{\theta}$. 
\end{theorem}

The next result on asymptotic normality of the adaptive MLEs requires condition (SD)$(\thb)$. 

\begin{theorem}\label{theo5-2}
Assume {\rm(a1)}, {\rm(a2')}, $({\rm GLM}^*)$ with (\ref{eq-a3'}), and $({\rm SD})(\thb)$. Assume further that the inverse link function $G$ 
is twice continuously differentiable, $\thb\in{\rm int}(\Theta)$, and the adaptive estimators employed by the algorithm are strongly consistent, that is,
$\widehat{\theta}_k\asto\thb$. Then 
\[
\sqrt{n_k}\,\bigl(\widehat{\theta}_k^{(\rm\scriptsize ML)}-\thb\bigr)\dto {\rm N}\bigl(0,\Ms^{-1}(\thb)\bigr), 
\]
where $\Ms(\thb)$ is given by condition $({\rm SD})(\thb)$.
\end{theorem}

\noindent
{\em Example 5: Simulation.}\\
We illustrate the results on consistency and asymptotic normality of the maximum likelihood estimators 
(Theorems \ref{theo5-1} and \ref{theo5-2}) and the asymptotic D-optimality of the generated designs (Corollary \ref{cor3-1})  
by simulations under the logit model (i) of Example 3 in Section 2. 
The experimental interval was chosen as ${\cal X}=[-4\,,\,4\,]$  and the parameter space  
as a rectangle $\Theta=[-10\,,\,10\,]\times [\,0.1\,,\,10\,]$. By simulations   
$10{,}000$ paths (more precisely: pieces of paths up to $k=250$) of the 2-step algorithm were generated 
for each of two cases of true parameter points: $\thb=(0,1)^\trp$ and $\thb=(4,1)^\trp$. The maximum likelihood estimators 
were employed, that is, $\widehat{\theta}_k=\MLk$. The starting design $\xi_1$ was always the three point design
with support points $-4$, $0$, $4$ and uniform weights $1/3$. So after step $k$ the total number of 
observations included is $n=n_k=2k+1$. In fact, $n$ rather than $k$ is used when comparing 
to the adaptive Wynn algorithm which is a 1-step algorithm. To this end, also
$10{,}000$ paths of the adaptive Wynn algorithm employing adaptive maximum likelihood estimates were simulated,
again for each of the two cases $\thb=(0,1)^\trp$ and $\thb=(4,1)^\trp$.        
Addressing the (almost sure) asymptotic D-optimality of the designs generated by the algorithms 
the development (as $n$ grows) of the D-efficiencies of the generated designs from the simulated paths 
is focussed (see top pictures in Figure 1). The D-efficiency (at the true parameter point $\thb$) of a design $\xi$ is defined by
$\bigl\{\det M(\xi,\thb)\big/\det M_*(\thb)\bigr\}^{1/2}$, where $M_*(\thb)$ is the information matrix 
of the locally D-optimal design at $\thb$. For the two cases of $\thb$ considered here the locally D-optimal design
$\xi^*_{\scriptsize\thb}$ and the inverse of its information matrix at $\thb$ are given by 
\begin{eqnarray*}
&&\thb=(0,1)^\trp\,:\ \ \xi^*_{\scriptsize\thb}=\left({\small\begin{array}{cc}-1.543 & 1.543 \\ 1/2 & 1/2\end{array}}\right),      
\quad M_*^{-1}(\thb) = \left({\small \begin{array}{cc} 6.899 & 0 \\ 0 & 2.894\end{array}}\right);\\
&&\thb=(4,1)^\trp\,:\ \ \xi^*_{\scriptsize\thb}=
\left({\small\begin{array}{cc}-4 & -1.601 \\ 1/2 & 1/2\end{array}}\right),      
\quad M_*^{-1}(\thb) = \left({\small \begin{array}{cc} 76.415 & 20.438 \\ 20.438 & 5.943\end{array}}\right),
\end{eqnarray*}
see Example 3 in Section 2 and Appendix A.1.
The consistency and asymptotic normality of the adaptive maximum likelihood estimators stated in Theorem \ref{theo5-1} and Theorem \ref{theo5-2}, 
respectively, should imply for the simulations that $n$ times the mean squared error matrix of the simulated parameter estimates
converges to $M_*^{-1}(\thb)$. This is illustrated in Figure 1 (middle plots) restricting to the diagonal entries of the matrices. 
Again, the adaptive 2-step {\em and} the adaptive Wynn algorithm are considered for a comparison. A further illustration of 
the asymptotic normality of maximum likelihood estimators from the adaptive 2-step algorithm is given by QQ-plots
in Figure 1 (bottom).               

The comparison of the two adaptive algorithms by our simulations suggests that both algorithms yield about the same convergence behavior
of the generated designs and maximum likelihood estimators. Note that the computation time of the adaptive Wynn 
was about double as large as that of the adaptive 2-step since the adaptive Wynn, as a `1-step ahead algorithm', carries out the optimization procedures (maximizing    
the likelihood function and the sensitivity)  twice as often as the adaptive 2-step.
However, for practical purposes it might be of greater importance that the adaptive 2-step algorithm allows some
parallel response sampling (two observations at a time) while the adaptive Wynn prescribes strictly sequential sampling  (one observation at   
a time). In particular, when observations are time consuming the adaptive 2-step may provide a substantial reduction 
of the total duration of data collection. 
\eop       

\begin{figure}[!tbp]
\vspace*{-0.045\textheight}
 
   \begin{minipage}[b]{0.52\textwidth} 
    \includegraphics[width=\textwidth]{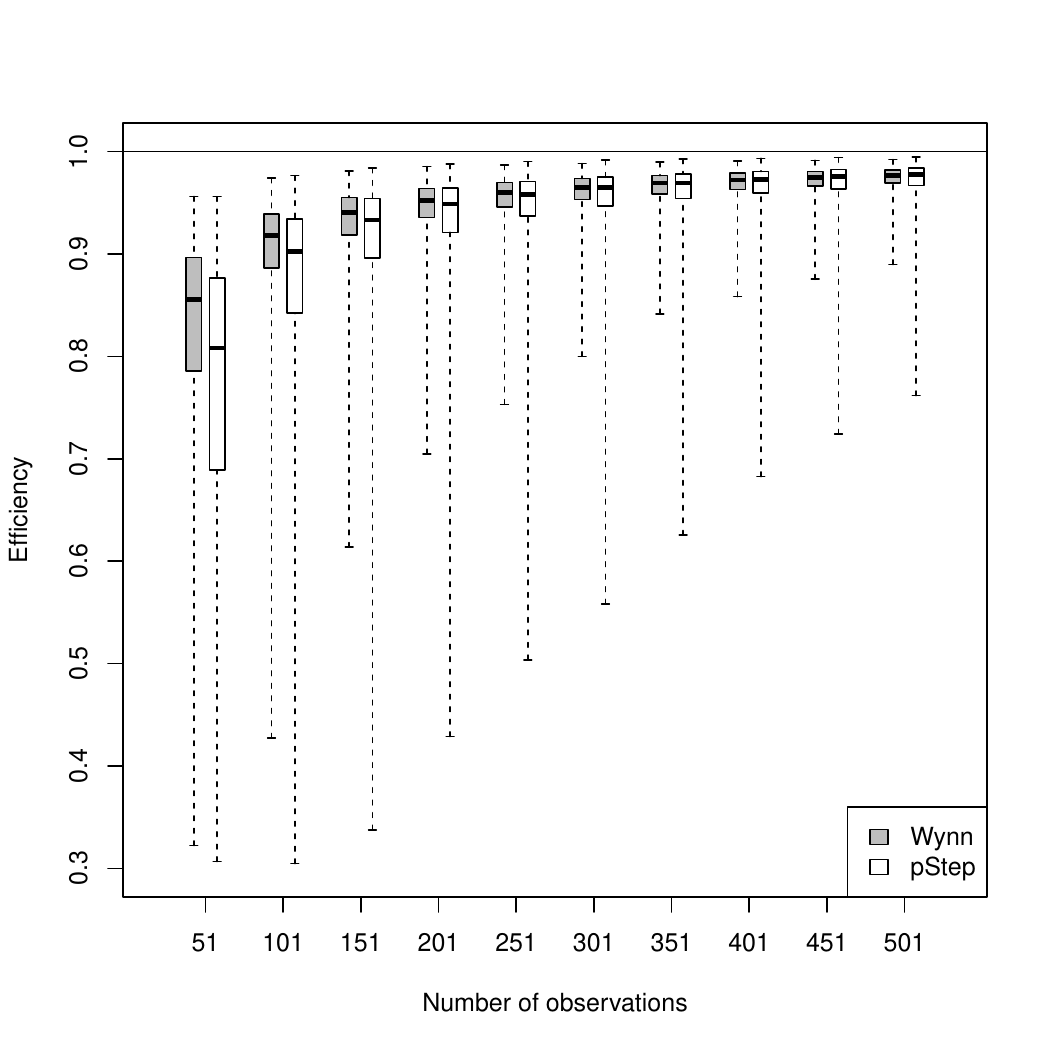}
   \end{minipage}
  \hfill 
  \begin{minipage}[b]{0.52\textwidth}
     \includegraphics[width=\textwidth]{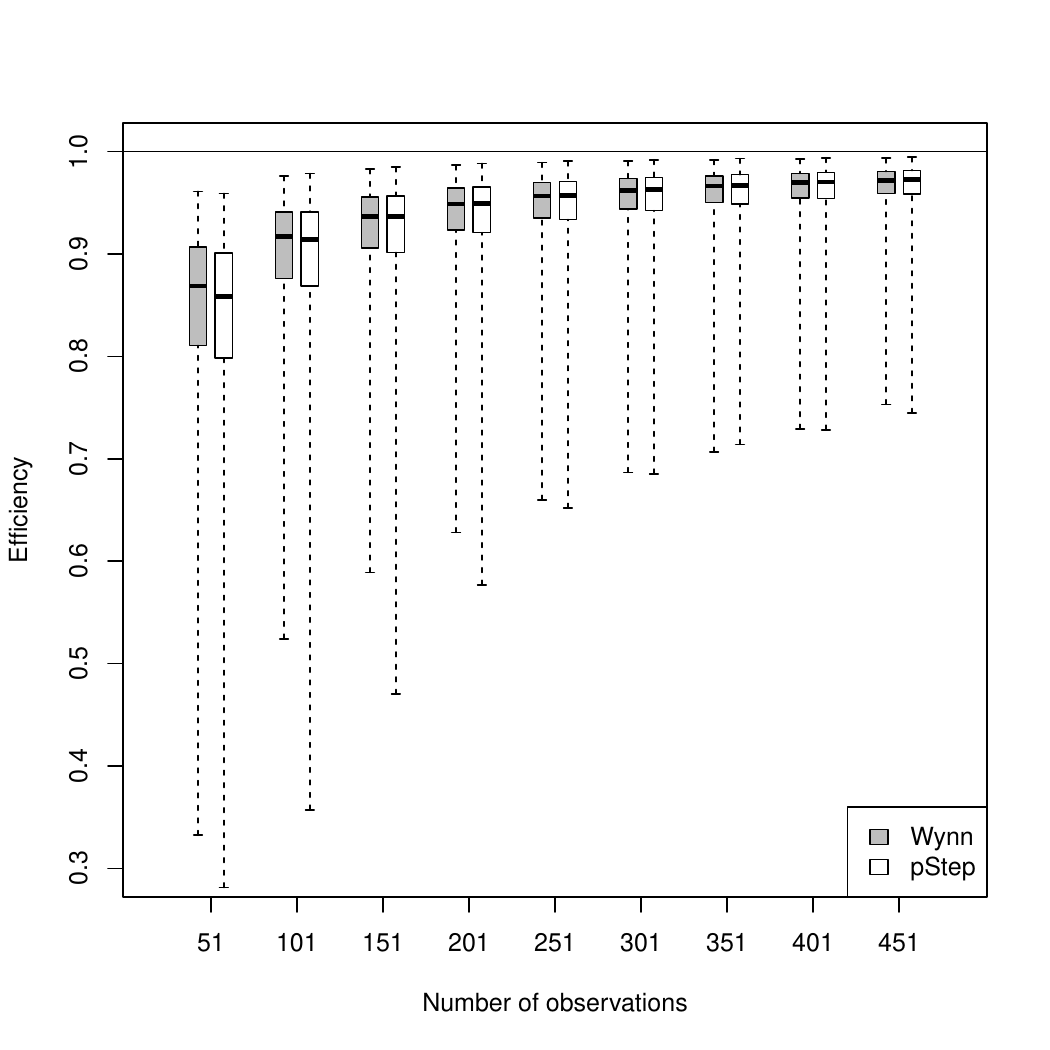}
  \end{minipage}
  
  \vspace*{-0.02\textheight}

  \begin{minipage}[b]{0.52\textwidth}
    \includegraphics[width=\textwidth]{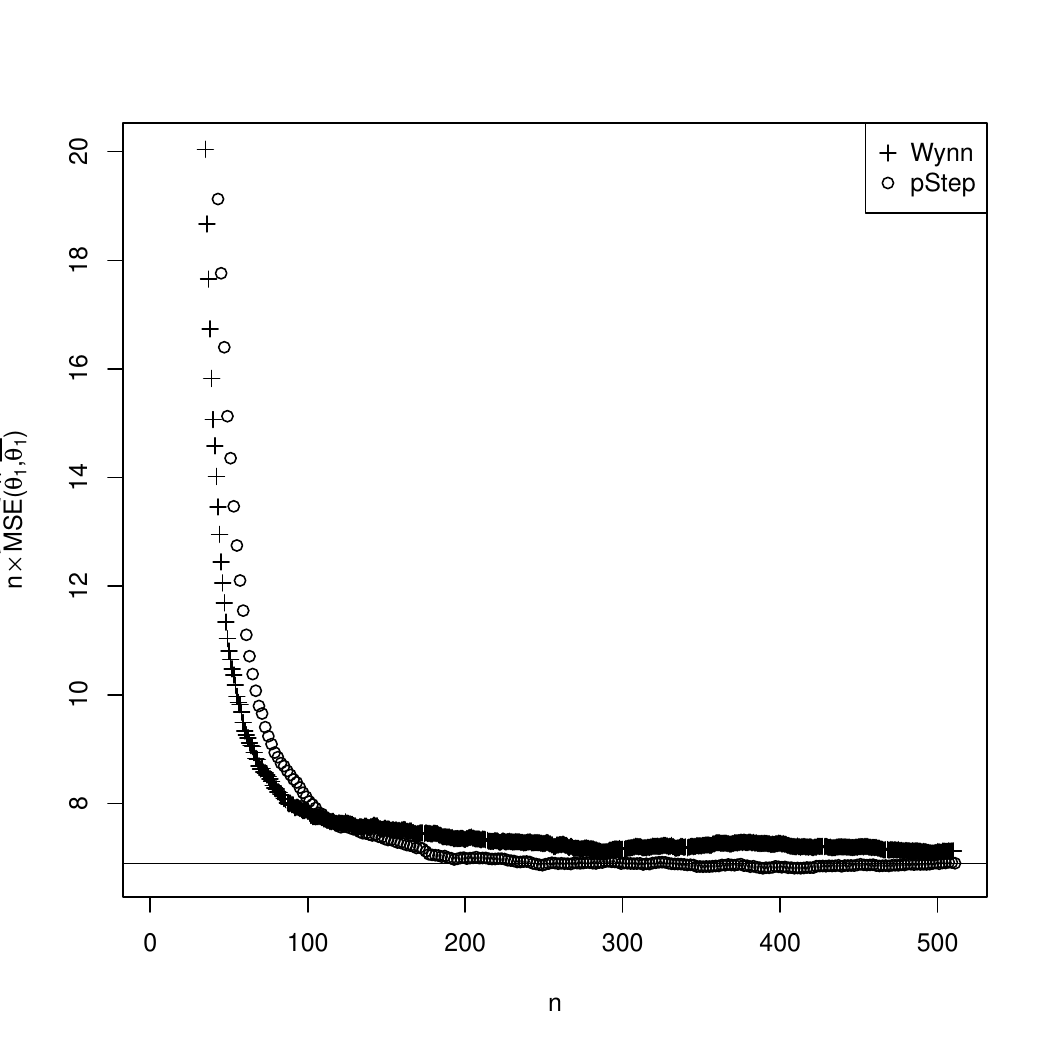}
  \end{minipage}
  \hfill
  \begin{minipage}[b]{0.52\textwidth}
    \includegraphics[width=\textwidth]{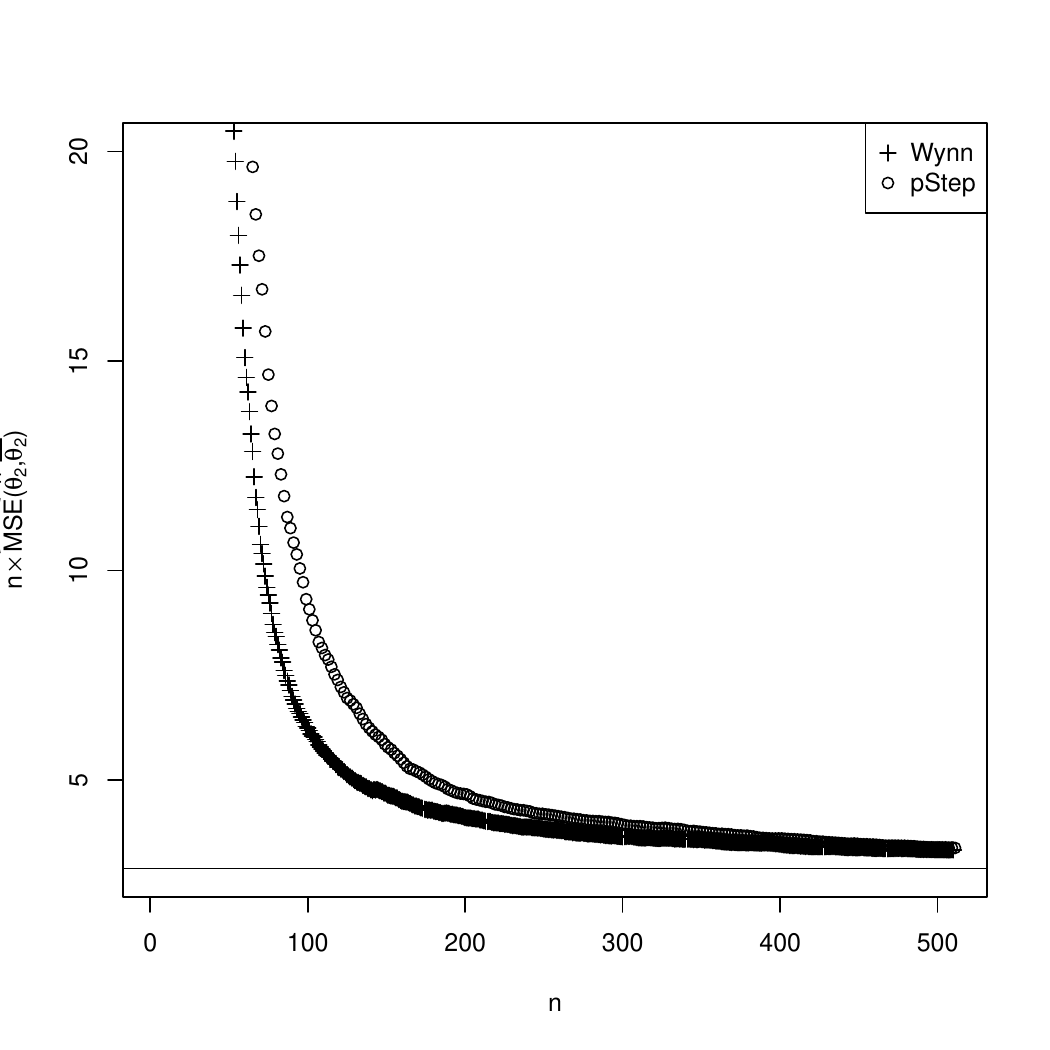}
  \end{minipage}

  \vspace*{-0.03\textheight}

  \begin{minipage}[b]{0.52\textwidth}
    \includegraphics[width=\textwidth]{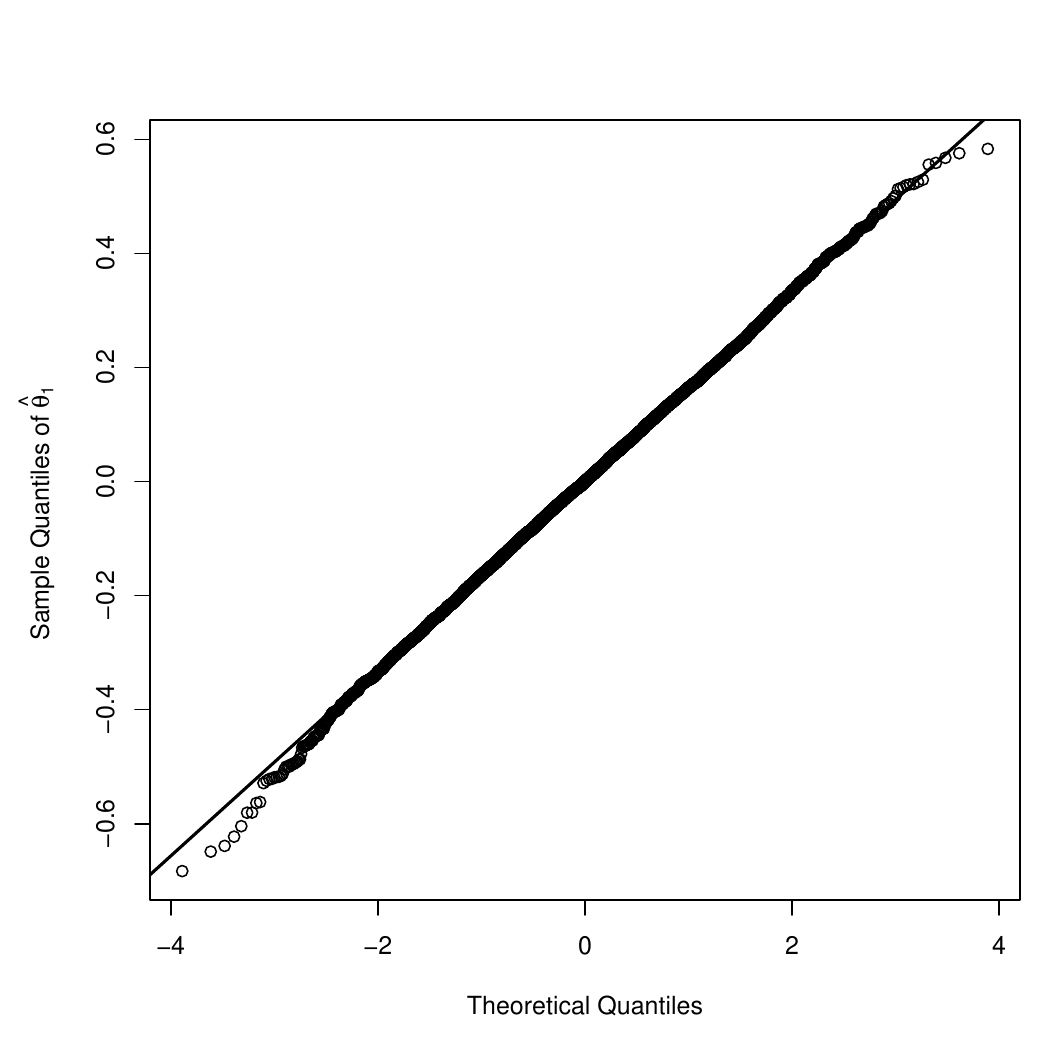}
  \end{minipage}
  \hfill
  \begin{minipage}[b]{0.52\textwidth}
    \includegraphics[width=\textwidth]{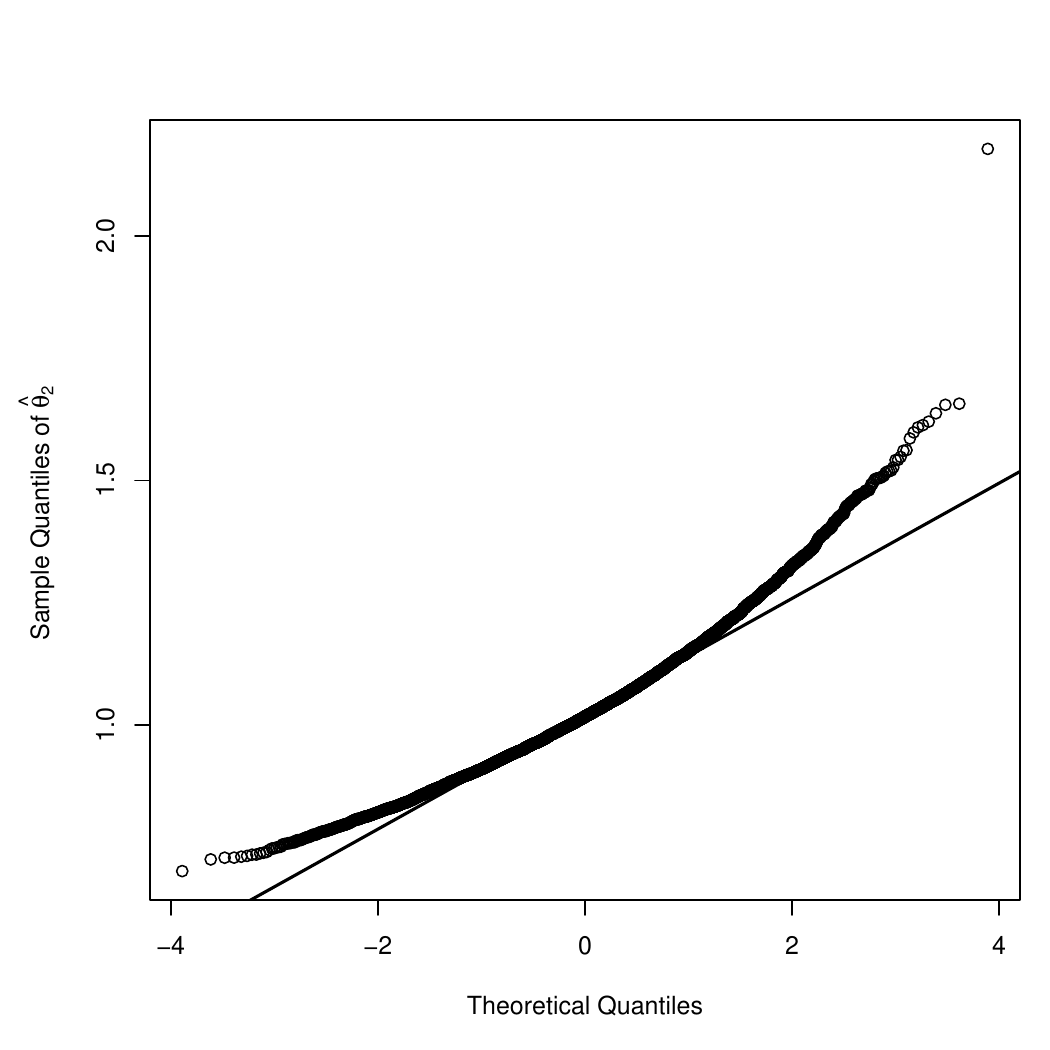}
  \end{minipage}

\vspace*{-0.015\textheight}
\caption{\small \underline{Top:} Box plots of D-efficiencies,  
where $\thb=(0,1)^\trp$ (left) and $\thb=(4,1)^\trp$ (right). The whiskers are from minimum to maximum.
\underline{Middle:} plots of $n$ times the mean squared error of the estimates of the first component (left) and  
the second component  (right)  of $\theta$, where $\thb=(0,1)^\trp$; the horizontal lines refer to the asymptotic variances (the diagonal entries
of $M_*^{-1}(\thb)$). 
\underline{Bottom:} adaptive 2-step algorithm and $\thb=(0,1)^\trp$; QQ-plots for the estimates of the components of 
$\theta$ at $k=125$, first component left, second component right.} 

\end{figure}

\begin{appendix}
\section{Appendix}
\setcounter{equation}{0}
\subsection{Supplement to Example 3.} 
Consider a model from Example 3 with a transformed design variable $z=\vartheta_1+\vartheta_2x$, where 
$\theta=(\vartheta_1,\vartheta_2)^\trp\in\Theta$ with $\vartheta_2\not=0$ is a given parameter point. Hence $z\in{\cal Z}=[\,\alpha\,,\,\beta\,]$,
the transform of  the design interval ${\cal X}=[\,a\,,\,b\,]$, that is, $\alpha<\beta$ are given by 
$\vartheta_1+\vartheta_2a$ and $\vartheta_1+\vartheta_2b$ arranged in increasing order.  
Denote $f_0(z)=\varphi(z)\,(1,z)^\trp$. Since, for any $z_1,z_2\in{\cal Z}$,
\[
\bigl\{\det\bigl[f_0(z_1)\,,\,f_0(z_2)\bigr]\bigr\}^2 = \varphi^2(z_1)\,\varphi^2(z_2)\,(z_2-z_1)^2
\]
a D-optimal saturated design (in the transformed model) can be viewed as an optimal solution to the problem 
\begin{eqnarray}
&&\mbox{maximize }\  \ h(z)=\ln \varphi(z_1) + \ln\varphi(z_2) + \ln(z_2-z_1)\ \mbox{ over }\ z\in D_{\alpha,\beta},\label{eqA-01}\\
&&\mbox{where } \ D_{\alpha,\beta}=[\,\alpha\,,\,\beta\,]^2\cap D,\ \ D=\bigl\{z=(z_1,z_2)\in\mathbb{R}^2:\,z_1<z_2\bigr\}.\nonumber
\end{eqnarray}
We will also look at the unbounded problem of saturated D-optimality without the bounds $\alpha$ and $\beta$, 
\begin{equation}
\mbox{maximize } \  h(z) \ \mbox{ over }\ z\in D.\label{eqA-02}
\end{equation}
An optimal solution to (\ref{eqA-01}) exists, whereas an optimal solution to (\ref{eqA-02}) may or may not exist in general. 
We will restrict to the case that (\ref{eqA-02}) has an optimal solution, see Lemma \ref{lemA-01} below. 
This includes the models in Ford, Torsney and Wu \cite{Ford-Torsney-Wu} listed in Table 4 of their paper along with the optimal solutions to
the unrestricted problems (\ref{eqA-02}). In particular, the models (i) to (iv) of our present Example 3 are included. 
Lemma \ref{lemA-01} below gives a description of the D-optimal saturated designs, that is, the optimal solutions to  
(\ref{eqA-01}), under the assumption that $\varphi$ is continuously differentiable and $\ln\varphi$ is strictly concave on $\mathbb{R}$.
This covers models  (i), (ii), and (iv) of our Example 3. The result of the lemma is not new: Table 3 of [FTW-1992] presents 
a slightly more general result. Unfortunately, the proof in Section 6.6 of that paper is incomplete and somewhat ambiguous:
the authors assume in their case (c) on p.~582 that $\varphi^2(u)\,(u-z_1)^2$ is non-decreasing in $u\in[\,z_1,\,\beta\,]$ for any $z_1\ge\alpha$.
But this is not necessarily met if $\alpha\le z_1<z_1^{**}$ even under log-concavity of $\varphi$. Here $z^{**}=(z_1^{**},z_2^{**})$
denotes the optimal solution to the unbounded problem  (\ref{eqA-02}). Similarly, in their case (d) on p.~582 the authors assume
that $\varphi^2(u)\,(z_2-u)^2$ is non-increasing in $u\in[\,\alpha\,,\,z_2]$ for any $z_2\le\beta$, but this may fail if
$z_2^{**}<z_2\le\beta$ despite log-concavity of $\varphi$. For example, for the logistic model (i) $\varphi^2(u)=
\exp(u)\bigl/\bigl(1+\exp(u)\bigr)^2$ and for $\alpha=-5$, $\beta=1.3$ case (c) of [FTW-1992] occurs, but for $z_1=-4$   
one finds that $\varphi^2(u)\,(u-z_1)^2$ is decreasing in $u$ when $1\le u \le 1.3$. 
The next lemma restates the result along with a proof, where our slightly stronger assumption of {\em strict} concavity of $\ln \varphi$ 
ensures uniqueness of the optimal solutions to problems (\ref{eqA-01}) and (\ref{eqA-02}). We denote the partial derivatives 
of $h$ on $D$ by $h_j'$, $j=1,2$, that is,  
\begin{equation}
h_1'(z)=(\ln\varphi)'(z_1) - \frac{1}{z_2-z_1},\ \   h_2'(z)=(\ln\varphi)'(z_2) + \frac{1}{z_2-z_1},\ \ z=(z_1,z_2)\in D,\label{eqA-00}
\end{equation}
where  $(\ln \varphi)'=\varphi'/\varphi$ is the derivative of $\ln \varphi$.

\begin{lemma} \label{lemA-01}
Let $\varphi$ be a positive and continuously differentiable function on the real line and such that $\ln \varphi$   
is strictly concave. Assume that there exists an optimal solution $z^{**}=(z_1^{**},z_2^{**})$ to the unbounded problem (\ref{eqA-02}).
 Then:\\  
{\bf(a)} The optimal solution $z^{**}=(z_1^{**},z_2^{**})$ to problem (\ref{eqA-02}) is unique, and $z^{**}$ is the unique solution
to the equations $h_j'(z)=0$, $j=1,2$, $z\in D$.\\[.5ex]
{\bf(b)} The optimal solution $z^*=(z_1^*,z_2^*)$ to problem (\ref{eqA-01}) is unique, and $z^*$ is obtained as follows,
where four cases are distinguished.\\[.5ex]
{\rm(1)}    
Let $\alpha\le z_1^{**}$ and $\beta\ge z_2^{**}$. Then $z^*=z^{**}$.\\[.3ex]
{\rm(2)} 
Let $\alpha\le z_1^{**}$ and $\beta<z_2^{**}$. If $h_1'(\alpha,\beta)\le0$ then $z^*=(\alpha,\beta)$; if
$h_1'(\alpha,\beta)>0$ then $z^*=(u,\beta)$ with $\alpha< u<\beta$ and $h_1'(u,\beta)=0$.\\[.3ex]  
{\rm(3)} 
 Let $\alpha> z_1^{**}$ and $\beta\ge z_2^{**}$. If $h_2'(\alpha,\beta)\ge0$ then $z^*=(\alpha,\beta)$; if
$h_2'(\alpha,\beta)<0$ then $z^*=(\alpha,v)$ with $\alpha< v<\beta$ and $h_2'(\alpha,v)=0$.\\[.3ex]  
{\rm(4)} 
Let $\alpha> z_1^{**}$ and $\beta< z_2^{**}$. Then $z^*=(\alpha,\beta)$.
\end{lemma} 

\noindent
{\bf Proof.} \ By strict concavity of $\ln \varphi$ the function $h$ is strictly concave on the convex set $D$.
Hence the optimal solution $z^{**}\in D$ to (\ref{eqA-02}) is unique and the optimal solution $z^*\in D_{\alpha,\beta}$ to (\ref{eqA-01}) is unique. 
Since $D$ is an open convex set, $z^{**}$ is the unique point of $D$ at which the gradient of $h$
is equal to zero, that is, $h_j'(z^{**})=0$, $j=1,2$. This proves part (a) of the lemma and uniqueness of $z^*$ in part (b). Case (1) of part (b)    
means that $z^{**}\in D_{\alpha,\beta}$ and hence $z^*=z^{**}$. 
In each of the remaining cases (2), (3), and (4) the optimal solution $z^*$ must have
at least one component equal to an end point of the interval $[\,\alpha\,,\,\beta\,]$, since otherwise $z^*$ would be an interior point of 
$D_{\alpha,\beta}$ entailing that the gradient of $h$ at $z^*$ equals zero and thus $z^*=z^{**}$. This is excluded in each of the cases (2), (3), and (4).
So, either $z^*=(u,\beta)$ with some $\alpha< u<\beta$, or $z^*=(\alpha,v)$ with some $\alpha< v <\beta$, or $z^*=(\alpha,\beta)$. 
As it is well-known, if $C\subseteq D$ is a given nonempty convex subset of $D$ then a point $\overline{z}=(\overline{z}_1,
\overline{z}_2)\in C$ maximizes $h(z)$ over $z\in C$ if and only if  the directional derivatives
of $h$ at $\overline{z}$ are nonpositive for all feasible directions, that is
\begin{equation}
h_1'(\overline{z})\,(z_1-\overline{z}_1) + h_2'(\overline{z})\,(z_2-\overline{z}_2)\,\le0\ \ \mbox{ for all }\ z=(z_1,z_2)\in C.
\label{eqA-03}
\end{equation}   
Consider case (2). Suppose that $z^*=(\alpha,v)$ with $\alpha<v<\beta$. From condition (\ref{eqA-03}) with $C=D_{\alpha,\beta}$ 
and $\overline{z}=z^*$ one gets
$h_1'(z^*)\le0$ and $h_2'(z^*)=0$. By (\ref{eqA-03}) with $C=D_\alpha=[\,\alpha\,,\,\infty)^2 \cap D$ and $\overline{z}=z^*$, one gets that $z^*$
maximizes $h(z)$ over $z\in D_\alpha$. But $z^{**}\in D_\alpha$ and thus $z^{**}$ is the unique maximizer of $h(z)$ over $z\in D_\alpha$.
Hence $z^*=z^{**}$ which is a contradiction. So the second component of $z^*$ must be equal to $\beta$ and, clearly, the first
component $z_1^*$ maximizes the function $z_1\mapsto h(z_1,\beta)$ over $z_1\in[\,\alpha\,,\,\beta\,)$. 
The derivative of that function is given by $h_1'(z_1,\beta)$, which is decreasing in $z_1\in[\,\alpha\,,\,\beta\,)$
and $h_1'(z_1,\beta)\to-\infty$ as $z_1\to\beta$.
One concludes: if $h_1'(\alpha,\beta)\le0$ then $z_1^*=\alpha$; otherwise $z_1^*=u$ the unique solution to $h_1'(u,\beta)=0$.\\
In case (3) the proof is analogous. Consider case (4). From (\ref{eqA-00}) it is obvious that $h_1'(z_1,z_2)$ is decreasing in
$z_1<z_2$ for fixed $z_2$, and increasing in $z_2>z_1$ for fixed $z_1$. Hence
\[
0=h_1'(z_1^{**},z_2^{**}) \ge h_1'(\alpha,z_2^{**})\ge h_1'(\alpha,\beta).
\]
Similarly, by (\ref{eqA-00}) the partial derivative $h_2'(z_1,z_2)$ is increasing in $z_1<z_2$ for fixed $z_2$, and decreasing in 
$z_2>z_1$ for fixed $z_1$. Hence
\[
0=h_2'(z_1^{**},z_2^{**}) \le h_2'(\alpha,z_2^{**})\le h_2'(\alpha,\beta).
\]
We have thus obtained that $h_1'(\alpha,\beta)\le0$ and $h_2'(\alpha,\beta)\ge0$. By condition (\ref{eqA-03})
with $C=D_{\alpha,\beta}$ and $\overline{z}=(\alpha,\beta)$ it follows that $z^*=(\alpha,\beta)$.
\eop

\subsection{Supplement to Example 4.}
Transforming the design variable and the design space, for a given parameter point $\theta=(\vartheta_1,\vartheta_2)^\trp$, 
to $z=(z_1,z_2)\in {\cal Z}=[\,0\,,\,c_1]\times[\,0\,,\,c_2]$ as described in Example 4, we have 
$f_0(z)=\overline{\varphi}_l(z)\,(1,z_1,z_2)^\trp$, where the index $l=1,2,3$ refers to the different cases (i), (ii), (iii) and
\begin{equation}
\overline{\varphi}_1(z)=\exp\bigl\{{\textstyle-\frac{1}{2}}(z_1+z_2)\bigr\},\ \ 
\overline{\varphi}_2(z)=\exp\bigl\{{\textstyle-\frac{1}{2}}z_1\bigr\},\ \
\overline{\varphi}_3(z)=1.\label{eqA-04}
\end{equation}
A D-optimal saturated design in the transformed model is described by three points $x^*=(x_1^*,x_2^*)$, $y^*=(y_1^*,y_2^*)$, and 
$z^*=(z_1^*,z_2^*)$ in the rectangle ${\cal Z}$ which maximize the function
\begin{eqnarray}
&&g_l(x,y,z) = \overline{\varphi}_l(x)\,\overline{\varphi}_l(y)\,\overline{\varphi}_l(z)\,\big|\det C(x,y,z)\big|,\label{eqA-05}\\ 
&&\mbox{where }\ \ C(x,y,z)=\,\left[\begin{array}{ccc}1 & 1 & 1\\x_1 & y_1 & z_1\\x_2& y_2 & z_2\end{array}\right]\nonumber
\end{eqnarray}
over all $x=(x_1,x_2),\, y=(y_1,y_2),\,z=(z_1,z_2)\,\in{\cal Z}$. 
Note that, again, the index $l=1,2,3$ refers to the different cases (i), (ii), and (iii). 
The next lemma shows the D-optimal saturated designs for each of the three cases.

\begin{lemma}\label{lemA-1}
Consider the functions $g_1$, $g_2$, and $g_3$ on ${\cal Z}^3$ defined  by (\ref{eqA-05}) and (\ref{eqA-04}), where ${\cal Z}=
[\,0\,,\,c_1]\times[\,0\,,\,c_2]$ with given $c_1>0$ and $c_2>0$. Then:\\[.5ex]
{\bf(i)} The points  which maximize $g_1(x,y,z)$ over 
$(x,y,z)\in{\cal Z}^3$ are the triples with components $(0,0)$, $(c_1^*,0)$, $(0,c_2^*)$ (arranged in any order), where
$c_j^*:=\min\{c_j,2\}$, $j=1,2$.\\[.5ex]
{\bf(ii)} The points  which maximize $g_2(x,y,z)$ over 
$(x,y,z)\in{\cal Z}^3$ are the triples with components $(0,0)$, $(c_1^*,\beta)$, $(0,c_2)$, where
$0\le\beta\le c_2$ is arbitrary and, as above, $c_1^*:=\min\{c_1,2\}$.\\[.5ex]
{\bf(iii)} The points  which maximize $g_3(x,y,z)$ over $\bigl(x,y,z)\in{\cal Z}^3$ 
are the triples with components $(0,0)$, $(c_1,0)$, $(\alpha,c_2)$,
the triples with components $(0,0)$, $(c_1,\beta)$, $(0,c_2)$,
the triples with components $(0,\beta)$, $(c_1,0)$, $(c_1,c_2)$,
and the triples with components $(\alpha,0)$, $(0,c_2)$, $(c_1,c_2)$,
where $0\le\alpha\le c_1$ and $0\le\beta\le c_2$ are arbitrary.
\end{lemma}

\noindent{\em Remark.} Geometrically, the solutions in part (iii) of the lemma are the triples consisting
of two adjacent vertices of the rectangle ${\cal Z}$ and any point from the edge of ${\cal Z}$ opposite
to the edge joining the two vertices.
\eop

\noindent
{\bf Proof.} \ Clearly, for $(x,y,z)\in{\cal Z}^3$ the product 
$\overline{\varphi}_l(x)\,\overline{\varphi}_l(y)\,\overline{\varphi}_l(z)$ in (\ref{eqA-05}) is equal to\\ 
$\exp\bigl\{{\textstyle-\frac{1}{2}\sum\limits_{j=1}^2(x_j+y_j+z_j)}\bigr\}$ in case $l=1$,
equal to $\exp\bigl\{{\textstyle-\frac{1}{2}(x_1+y_1+z_1)}\bigr\}$ in case  $l=2$, and equal to $1$ in case  $l=3$.
For later use, we show the following.
\begin{eqnarray}
&&\mbox{If $x,y,z\in{\cal Z}$ such that }  \min\{x_j,y_j,z_j\}=0 \ \mbox{ for $j=1,2$, \ then}\nonumber\\
&& \big|\det C(x,y,z)\big|\,\le\,\max\{x_1,y_1,z_1\}\,\max\{x_2,y_2,z_2\}.\label{eqA-06}
\end{eqnarray}
To see this, after a suitable permutation of $x$, $y$, and $z$, the following two cases have to be considered.
\underbar{Case 1:} \ $x_1=0$ and $y_2=0$;\quad\underbar{Case 2:} \ $x_1=0$ and $x_2=0$.\\
Assume Case 1. Then $\det C(x,y,z)=x_2(z_1-y_1)+y_1z_2$. If $y_1\le z_1$ then\\
$
\big|x_2(z_1-y_1)+y_1z_2\big| = x_2(z_1-y_1)+y_1z_2\le \max\{x_2,z_2\}\,(z_1-y_1+y_1)
= z_1\max\{x_2,z_2\}$, hence (\ref{eqA-06}). If $y_1> z_1$ then     
$
\big|x_2(z_1-y_1)+y_1z_2\big|\,\le\,\max\bigl\{x_2(y_1-z_1)\,,\,y_1z_2\bigr\}\,\le\,
y_1\,\max\{x_2,z_2\},
$
and hence (\ref{eqA-06}). Now assume Case 2. Then 
$
\big|\det C(x,y,z) \big| = \big|y_1z_2-z_1y_2\big|\,\le\,\max\{y_1z_2\,,\,z_1y_2\}
\le \max\{y_1,z_1\}\,\max\{y_2,z_2\},
$
and hence (\ref{eqA-06}).
Below we will use the fact that the function $\exp\bigl\{-\frac{1}{2}\,t\bigr\}\,t$ increases for $0\le t\le 2$ and decreases
for $2\le t<\infty$, and hence for $j=1,2$,
\begin{equation}
\exp\bigl\{-{\textstyle \frac{1}{2}}t\bigr\}\,t\,\le\,\exp\bigl\{-{\textstyle \frac{1}{2}}c_j^*\bigr\}\,c_j^*
\quad\mbox{for all $0\le t\le c_j$,}\label{eqA-07}
\end{equation} 
and the inequality is strict unless $t=c_j^*$.\\[.5ex] 
{\bf(i)} 
Let any $(x,y,z)\in{\cal Z}^3$ be given such that $C(x,y,z)$ is nonsingular. Denote $a_j=\min\{x_j,y_j,z_j\}$,
$j=1,2$, and $a=(a_1,a_2)$. Define $\widetilde{x}=x-a$, $\widetilde{y}=y-a$, and 
$\widetilde{z}=z-a$. Then $\widetilde{x},\widetilde{y},\widetilde{z}\in{\cal Z}$, 
$\det C(\widetilde{x},\widetilde{y},\widetilde{z})=\det C(x,y,z)$, and
$g_1(\widetilde{x},\widetilde{y},\widetilde{z})\ge g_1(x,y,z)$ with strict inequality unless $a_1=a_2=0$.
So, the maximizers of $g_1(x,y,z)$ over ${\cal Z}^3$ 
are among those points $(x,y,z)\in{\cal Z}^3$ such that $\min\{x_1,y_1,z_1\}=0$ and $\min\{x_2,y_2,z_2\}=0$.
Denote $g^*\,=\,\exp\bigl\{-\frac{1}{2}(c_1^*+c_2^*)\bigr\}\,c_1^*c_2^*$ which is the common value $g_1(x^*,y^*,z^*)$
of the claimed maximizers $(x^*,y^*,z^*)$ in part (i) of the lemma. Let $(x,y,z)\in{\cal Z}$ be given 
such that $\min\{x_j,y_j,z_j\}=0$ for $j=1,2$. Then, by (\ref{eqA-06}),
and (\ref{eqA-07}), 
\begin{eqnarray*}
g_1(x,y,z) &\le & \exp\bigl\{-{\textstyle\frac{1}{2}}\bigl(\max\{x_1,y_1,z_1\}+\max\{x_2,y_2,z_2\}\bigr)\bigr\}\,   
\max\{x_1,y_1,z_1\}\,\max\{x_2,y_2,z_2\}\\
&\le& \exp\bigl\{-{\textstyle\frac{1}{2}}(c_1^*+c_2^*)\bigr\}\,c_1^*c_2^*\,=\,g^*,
\end{eqnarray*}
and the equality $g_1(x,y,z)=g^*$ implies that for each $j=1,2$ two of $x_j,y_j,z_j$ must be zero, 
from which one concludes $\{x,y,z\}=\bigl\{(0,0)\,,\,(c_1^*,0)\,,\,(0,c_2^*)\bigr\}$.\\[.5ex]
{\bf(ii)} Let any $(x,y,z)\in{\cal Z}^3$ be given such that $C(x,y,z)$ is nonsingular. For $j=1,2$ denote\\ 
$a_j=\min\{x_j,y_j,z_j\}$, and $\overline{a}_2=\max\{x_2,y_2,z_2\}$ and $\lambda=c_2/(\overline{a}_2-a_2)$. 
Note that $a_2<\overline{a}_2\le c_2$ and $\lambda\ge1$. Define 
$\widetilde{x}=\bigl(x_1-a_1,\lambda(x_2-a_2)\bigr)$, $\widetilde{y}=\bigl(y_1-a_1,\lambda(y_2-a_2)\bigr)$,
and $\widetilde{z}=\bigl(z_1-a_1,\lambda(z_2-a_2)\bigr)$. Then,  $(\widetilde{x},\widetilde{y},\widetilde{z})\in{\cal Z}^3$
and $\det C(\widetilde{x},\widetilde{y},\widetilde{z}) = \lambda\det C(x,y,z)$ and $\widetilde{x}_1+\widetilde{y_1}+\widetilde{z}
\le x_1+y_1+z_1$. It follows that  $g_2(\widetilde{x},\widetilde{y},\widetilde{z})\ge g_2(x,y,z)$ and the inequality is strict unless
$a_1=a_2=0$ and $\overline{a}_2=c_2$. So the maximizers of $g_2(x,y,z)$ over $(x,y,z)\in{\cal Z}^3$ must satisfy
$\min\{x_j,y_j,z_j\}=0$ for $j=1,2$, and $\max\{x_2,y_2,z_2\}=c_2$. Let $(x,y,z)\in{\cal Z}^3$ be 
any such triple. By (\ref{eqA-06}) and (\ref{eqA-07})
\begin{eqnarray*}
g_2(x,y,z) &\le & \exp\bigl\{-{\textstyle\frac{1}{2}}\max\{x_1,y_1,z_1\}\bigr)\bigr\}\,   
\max\{x_1,y_1,z_1\}\,c_2\\
&\le& \exp\bigl\{-{\textstyle\frac{1}{2}}c_1^*\bigr\}\,c_1^*c_2\,=\,g^*,
\end{eqnarray*}
and the equality $g_2(x,y,z)=g^*$ implies that two of $x_1,y_1,z_1$ are equal to zero and one is  equal to $c_1^*$,
and $\max\{x_2,y_2,z_2\}=c_2$, $\min\{x_2,y_2,z_2\}=0$. From this the result of part (ii) follows.\\[.5ex]
{\bf(iii)} Let any $(x,y,z)\in{\cal Z}^3$ be given such that $C(x,y,z)$ is nonsingular. With $a_j=\min\{x_j,y_j,z_j\}$,
$\overline{a}_j=\max\{x_j,y_j,z_j\}$ and $\lambda_j=c_j/(\overline{a}_j-a_j)$, $j=1,2$, define
$\widetilde{x}=\bigl(\lambda_1(x_1-a_1),\lambda_2(x_2-a_2)\bigr)$,
$\widetilde{y}=\bigl(\lambda_1(y_1-a_1),\lambda_2(y_2-a_2)\bigr)$, and
$\widetilde{z}=\bigl(\lambda_1(z_1-a_1),\lambda_2(z_2-a_2)\bigr)$.      
Then $\widetilde{x},\widetilde{y},\widetilde{z}\in{\cal Z}$, 
$\det C(\widetilde{x},\widetilde{y},\widetilde{z})=\lambda_1\lambda_2\det C(x,y,z)$, and since $\lambda_j\ge1$ for $j=1,2$
one has $g_3(\widetilde{x},\widetilde{y},\widetilde{z})\ge g_3(x,y,z)$ with strict inequality unless $a_j=0$ and 
$\overline{a}_j=c_j$ for $j=1,2$. So the maximizers of $g_3(x,y,z)$ are among the triples $(x,y,z)\in{\cal Z}^3$
such that after a suitable permutation of $x,y,z$, one has $(x_1,y_1,z_1)=(0,c_1,\alpha)$ 
for some $0\le\alpha\le c_1$, and $x_2,y_2,z_2$ is some (other) permutation of $0,c_2,\beta$
for some  $0\le\beta\le c_2$. Checking the six possible permutations and maximizing $g_3(x,y,z)$ w.r.t. the remaining variables
$\alpha$ and $\beta$ the four types of triples as stated in part (iii) of the lemma appear as the maximizers. 
\eop

In case (i) the lemma yields the uniformly weighted design on the three points $(0,0)$, $(c_1^*,0)$, $(0,c_2^*)$  
as the unique D-optimal saturated design in the transformed model $f_0(z)=\exp\bigl\{-\frac{1}{2}(z_1+z_2)\bigr\}\,(1,z_1,z_2)^\trp$,
$z=(z_1,z_2)\in{\cal Z}=[\,0\,,\,c_1]\times [\,0\,,\,c_2]$. Let us denote this design by $\zeta^*$. We ask whether $\zeta^*$ is D-optimal
(in the class of all designs $\zeta$). Some answers are given by the next lemma the first part of which is covered
by a more general result in \cite{Russell-et-al}, see the lemma on p.~723 of that paper. We present though an alternative (short) proof.

\begin{lemma}\label{lemA-3}
Assume case (i) and the transformed model, and let $\zeta^*$ be the design with support points  $(0,0)$, $(c_1^*,0)$, $(0,c_2^*)$
and uniform weights $1/3$. If $c_j\ge2$ for $j=1,2$ then $\zeta^*$ is D-optimal. If $c_j\le2$ for $j=1,2$ and
$\bigl(1+\exp(-c_1)\bigr)\,\bigl(1+\exp(-c_2)\bigr)\,>2$ then $\zeta^*$ is not D-optimal.
\end{lemma}   

\noindent
{\bf Proof.} \  Denoting
$C=\bigl[f_0(0,0)\,,\,f_0(c_1^*,0)\,,\,f_0(0,c_2^*)\bigr]$, the information matrix of $\zeta^*$ is given by $(1/3)CC^\trp$.
Denote $g_0(z)=C^{-1}f_0(z)$, $z\in{\cal Z}$. The condition for D-optimality of $\zeta^*$ from the Kiefer-Wolfowitz equivalence theorem
can be written as
\begin{equation}
g_0^\trp(z)\,g_0(z)\,\le 1\quad\mbox{for all $z\in{\cal Z}$.}\label{eqA-08}
\end{equation}
By straightforward calculation,\\ 
$g_0(z)=\exp\bigl(-(z_1+z_2)/2\bigr)\,\bigl(1-(z_1/c_1^*) -(z_2/c_2^*),\,\exp(c_1^*/2)\,z_1/c_1^*,\,\exp(c_2^*/2)\,z_2/c_2^*\bigr)^\trp$\\
for all $z=(z_1,z_2)\in{\cal Z}$, and hence
\begin{equation}
g_0^\trp(z)\,g_0(z) = 
\exp\bigl(-(z_1+z_2)\bigr)\,\Bigl\{\Bigl(1-\frac{z_1}{c_1^*} - \frac{z_2}{c_2^*}\Bigr)^2 + \exp(c_1^*)\,\Bigl(\frac{z_1}{c_1^*}\Bigr)^2 
+ \exp(c_2^*)\,\Bigl(\frac{z_2}{c_2^*}\Bigr)^2\Bigr\}. 
\label{eqA-09}
\end{equation}
Consider the case that $c_j\ge2$ for $j=1,2$, hence $c_j^*=2$ for $j=1,2$. To prove (\ref{eqA-08}), observing (\ref{eqA-09}), 
it suffices to show that for any $t\in[\,0\,,\,\infty]$ one has 
\begin{eqnarray}
&&e^{-t}\,\Bigl\{\Bigl(1-{\textstyle\frac{1}{2}}z_1- {\textstyle\frac{1}{2}}z_2\Bigr)^2 +   
{\textstyle\frac{e^2}{4}}\,\bigl(z_1^2 + z_2^2\bigr)\Bigr\} \ \le 1\ \mbox{ for all $z\in H_t$,}\label{eqA-010}\\ 
&&\mbox{where }\ H_t=\bigl\{z=(z_1,z_2)\in[\,0\,,\,\infty)^2\,:\,z_1+z_2=t\bigr\}.\nonumber 
\end{eqnarray}
On $H_t$, for any fixed $t$, the function on the l.h.s. of (\ref{eqA-010}) is convex and hence attains its maximum at
an end-point of the line segment $H_t$, that is, at $(t,0)$ or $(0,t)$ whose common function value is
\[
h(t) = e^{-t}\,\Bigl\{\Bigl(1-{\textstyle\frac{1}{2}}t\Bigr)^2 +   
{\textstyle\frac{e^2}{4}}\,t^2\Bigr\}.
\]
We have to show that $h(t)\le1$ for all $t\ge0$. For the derivative of $h$ one calculates
\[
h'(t)= -{\textstyle\frac{1+e^2}{4}}\,e^{-t}\,\Bigl\{t^2 - 2\,{\textstyle\frac{3+e^2}{1+e^2}}\,t + {\textstyle\frac{8}{1+e^2}}\Bigr\}.
\]
It is easily seen that  $h'(t)\le0$ for $t\le t_1$, $h'(t)\ge0$ for $t_1\le t\le t_2$,  and $h'(t)\le0$ for $t\ge t_2$, 
where $t_1=4\big/\bigl(1+e^2\bigr)$ 
and $t_2=2$ are the zeros of $h'$. From this and by $h(0)=h(2)=1$, one gets $h(t)\le \max\{h(0),\,h(2)\}=1$ for all $t\ge0$.

Now consider the case that $c_j\le2$ for $j=1,2$, hence $c_j^*=c_j$ for $j=1,2$. For the vertex $\overline{z}=(c_1,c_2)$ of ${\cal Z}$
one gets by (\ref{eqA-09})
\[
g_0^\trp(\overline{z})\,g_0(\overline{z}) = 
\exp\bigl(-(c_1+c_2)\bigr)\,\bigl\{1 + \exp(c_1) + \exp(c_2)\bigr\} = \bigl(1+\exp(-c_1)\bigr)\,\bigl(1+\exp(-c_2)\bigr)\,-1.
\]
So, if $\bigl(1+\exp(-c_1)\bigr)\,\bigl(1+\exp(-c_2)\bigr)\,>2$ then the Kiefer-Wolfowitz condition (\ref{eqA-08})
is violated at $z=\overline{z}$ and thus $\zeta^*$ is not D-optimal.
\eop

\subsection{Proofs of the lemmas}
{\bf Proof of Lemma \ref{lem3-1}}\\[.5ex] 
By the uniform continuity of the function $(x,\theta)\mapsto f_\theta(x)\,f_\theta^\trp(x)$
on the compact metric space ${\cal X}\times\Theta$ 
and by the assumption that $\lim_{k\to\infty}{\rm d}_\Theta(\rho_k,\tau_k)=0$, one gets 
\begin{equation}
\lim_{k\to\infty}\Bigl(\max_{x\in{\cal X}}\big\Vert f_{\rho_k}(x)\,f_{\rho_k}^\trp(x)\,-\,
f_{\tau_k}(x)\,f_{\tau_k}^\trp(x)\big\Vert\Bigr)\,=\,0.
\label{eq3-1}
\end{equation}
Let any $\xi\in\Xi$ be given. Then 
\begin{eqnarray*}
&&\big\Vert M(\xi,\rho_k)\,-\,M(\xi,\tau_k)\big\Vert\,=\,
\Big\Vert \sum_{x\in{\rm\scriptsize supp}(\xi)}\xi(x)\bigl[f_{\rho_k}(x)\,f_{\rho_k}^\trp(x)\,-\,
f_{\tau_k}(x)\,f_{\tau_k}^\trp(x)\bigr]\Big\Vert\\
&&\le\,\sum_{x\in{\rm\scriptsize supp}(\xi)}\xi(x)\,\big\Vert f_{\rho_k}(x)\,f_{\rho_k}^\trp(x)\,-\,
f_{\tau_k}(x)\,f_{\tau_k}^\trp(x)\big\Vert\\
&&\le\,
\max_{x\in{\rm\scriptsize supp}(\xi)}\big\Vert f_{\rho_k}(x)\,f_{\rho_k}^\trp(x)\,-\,
f_{\tau_k}(x)\,f_{\tau_k}^\trp(x)\big\Vert.
\end{eqnarray*}
Hence
\[
\sup_{\xi\in\Xi}\big\Vert M(\xi,\rho_k)\,-\,M(\xi,\tau_k)\big\Vert\,
\le\,\max_{x\in{\cal X}}\big\Vert f_{\rho_k}(x)\,f_{\rho_k}^\trp(x)\,-\,
f_{\tau_k}(x)\,f_{\tau_k}^\trp(x)\big\Vert,
\]
and together with (\ref{eq3-1}) the first result follows. We observe that
the set of all information matrices, 
${\cal M}=\bigl\{M(\xi,\theta)\,:\,\xi\in\Xi,\ \theta\in\Theta\bigl\}$,
is bounded since for all $\xi\in\Xi$ and $\theta\in\Theta$  
\[
\Vert M(\xi,\theta)\Vert\le
\max_{x\in{\cal X}}\Vert f_{\theta}(x)\,f^\trp_\theta(x)\Vert=\max_{x\in{\cal X}}\Vert f_\theta(x)\Vert^2\le \gamma^2
\]
where $\gamma:=\max_{(x,\beta)\in{\cal X}\times\Theta}\Vert f_\beta(x)\Vert<\infty$. 
So there is a compact subset ${\cal A}$ of the set of all nonnegative definite $p\times p$ matrices such that
${\cal M}\subseteq{\cal A}$. Now the second statement of the lemma follows using the uniform continuity of $\Phi$
on the compact set ${\cal A}$.  
\eop 

\newpage\noindent
{\bf Proof of Lemma \ref{lem3-2}}\\[.5ex] 
By (\ref{eq2-1}) and
$\eta_k=\frac{1}{p}\sum_{j=1}^p\delta[x_{n_k+j}]$, one gets
\begin{equation}
\xi_k\, =\, \frac{n_1}{n_k}\xi_1 \,+\,\frac{p}{n_k}\sum_{j=1}^{k-1} \eta_j,\ \ \mbox{for all $k\ge2$},\label{eq3-2}
\end{equation}
and hence, for all $\theta\in\Theta$ and all $k\ge2$, 
\begin{equation}
M(\xi_k,\theta)=\frac{n_1}{n_k}M(\xi_1,\theta) \,+\,\frac{p}{n_k}\sum_{j=1}^{k-1} M(\eta_j,\theta).
\label{eq3-3}
\end{equation}
Let $\eta^*$ be any saturated design which maximizes $\det M(\eta,\overline{\theta})$ over all saturated designs $\eta$.
Since $\eta_k$ maximizes $\det M(\eta,\theta_k)$ over all saturated designs $\eta$, one has 
\[
\det M(\eta_k,\theta_k) \ge \det M(\eta^*,\theta_k)\ 
\mbox{ for all $k\ge1$}.
\]
For $k\to\infty$ the r.h.s. of the latter inequality converges to $M(\eta^*,\overline{\theta})$ 
by Lemma \ref{lem3-1}, and hence
\[
\liminf_{k\to\infty} 
\det M(\eta_k,\theta_k) \ge \det M(\eta^*,\overline{\theta}).
\]
Again by Lemma \ref{lem3-1},   
\ $\det M(\eta_k,\theta_k) -\det M(\eta_k,\overline{\theta})\,\to\, 0$ \ as $k\to\infty$,
hence\\ 
$\liminf_{k\to\infty}\det M(\eta_k,\theta_k)$ and $\liminf_{k\to\infty}\det M(\eta_k,\overline{\theta})$
coincide, and thus
\[
\liminf_{k\to\infty}\det M(\eta_k,\overline{\theta})\,\ge\,    
\det M(\eta^*,\overline{\theta}).
\]
On the other hand, 
$\det M(\eta_k,\overline{\theta}) \le \det M(\eta^*,\overline{\theta})$ 
for all $k$ and hence 
\[
\limsup_{k\to\infty}\det M(\eta_k,\overline{\theta})\,\le\,    
\det M(\eta^*,\overline{\theta}).
\]
It follows that
\begin{equation}
\lim_{k\to\infty}\det M(\eta_k,\overline{\theta})\,=\,    
\det M(\eta^*,\overline{\theta}).\label{eq3-4}
\end{equation}
Denoting by $d_{{\rm\scriptsize s}*}(\overline{\theta})$ the common value of the determinants on 
${\cal M}_{{\rm\scriptsize s}*}(\overline{\theta})$, we have obtained that  
$\lim_{k\to\infty}\det M(\eta_k,\overline{\theta}) = d_{{\rm\scriptsize s}*}(\overline{\theta})$.
Consider the compact subset of information matrices
\[
{\cal M}_{{\rm\scriptsize s,u}}(\overline{\theta})\,=\,\Bigl\{\frac{1}{p}\sum_{j=1}^pf_{\thb}^\trp(z_j)\,f_{\thb}(z_j)\,:\,
z_1,\ldots,z_p\in{\cal X}\Bigr\},
\] 
which constitutes the closure of the set of information matrices at $\thb$ of all saturated designs with uniform weights.
For any given $\varepsilon>0$, consider the compact (or empty) subset of that set 
\[ 
\bigl\{M\in {\cal M}_{{\rm\scriptsize s,u}}(\overline{\theta})\,:\,
{\rm dist}\bigl(M,\,{\cal M}_{{\rm\scriptsize s}*}(\overline{\theta})\bigr)\ge\varepsilon\bigr\}.
\]
The maximum value of the determinant on the latter set (where $\max\emptyset:=-\infty$) is strictly less than
$d_{{\rm\scriptsize s}*}(\overline{\theta})$, and therefore  
${\rm dist}\bigl(M(\eta_k,\thb),\,{\cal M}_{{\rm\scriptsize s}*}(\overline{\theta})\bigr)<\varepsilon$ for $k$ large enough.
We have thus shown that
\begin{equation}
{\rm dist}\bigl(M(\eta_k,\overline{\theta}),\,{\cal M}_{{\rm\scriptsize s}*}(\overline{\theta})\bigr)\,\longrightarrow\,0
\ \mbox{ as $k\to\infty$.}
\label{eq3-4a}
\end{equation}
Trivially, (\ref{eq3-4a}) remains true when enlarging the set ${\cal M}_{{\rm\scriptsize s}*}(\overline{\theta})$
to its convex hull. Since the function $B\mapsto{\rm dist}\Bigl(B,\,
{\rm Conv}\,{\cal M}_{{\rm\scriptsize s}*}(\overline{\theta})\Bigr)$       
is convex on $\mathbb{R}^{p\times p}$, one gets for all $k\ge2$, 
\[
{\rm dist}\Bigl(\frac{1}{k-1}\sum_{j=1}^{k-1}M(\eta_j,\overline{\theta}),\,{\rm Conv}\,{\cal M}_{{\rm\scriptsize s}*}(\overline{\theta})
\Bigr) \,\le\,
\frac{1}{k-1}\sum_{j=1}^{k-1}{\rm dist}\Bigl(M(\eta_j,\overline{\theta}),\,{\rm Conv}\,{\cal M}_{{\rm\scriptsize s}*}(\overline{\theta})
\Bigr),
\]
and the r.h.s. goes to zero as $k\to\infty$ by (\ref{eq3-4a}). It follows that
\begin{equation}
{\rm dist}\Bigl(\frac{1}{k-1}\sum_{j=1}^{k-1}M(\eta_j,\overline{\theta}),
{\rm Conv}\,{\cal M}_{{\rm\scriptsize s}*}(\overline{\theta})\Bigr)\,\longrightarrow\,0\ \mbox{ as $k\to\infty$.}
\label{eq3-4b}
\end{equation}
By (\ref{eq3-3}), observing $n_k=n_1+(k-1)p$,
\[
\Big\Vert M(\xi_k,\thb)\,-\,\frac{1}{k-1}\sum_{j=1}^{k-1}M(\eta_j,\thb)\Big\Vert\,\longrightarrow0 \ \mbox{ as $k\to0$,}
\]
and together with (\ref{eq3-4b}), 
\begin{equation}
{\rm dist}\Bigl(M(\xi_k,\overline{\theta}),\,{\rm Conv}\,{\cal M}_{{\rm\scriptsize s}*}(\overline{\theta})\Bigr)\,\longrightarrow\,0
\ \mbox{ as $k\to\infty$.}\label{eq3-4c}
\end{equation}
Let $\theta_k'$, $k\in\mathbb{N}$, be any sequence in $\Theta$ which converges to $\overline{\theta}$.
Then, by (\ref{eq3-4c}) and Lemma \ref{lem3-1}, 
\[
{\rm dist}\Bigl(M(\xi_k,\theta_k'),\,{\rm Conv}\,{\cal M}_{{\rm\scriptsize s}*}(\overline{\theta})\Bigr)\,\longrightarrow\,0
\ \mbox{ as $k\to\infty$,}
\]
If condition (SD)$(\overline{\theta})$ holds then 
${\rm Conv}\,{\cal M}_{{\rm\scriptsize s}*}(\overline{\theta}) = \bigl\{M_{{\rm\scriptsize s}*}(\overline{\theta})\bigr\}$ and hence\\
$\lim_{k\to\infty}M(\xi_k,\theta_k')\,=\,M_{{\rm\scriptsize s}*}(\overline{\theta})$.
\eop

\noindent
{\bf Proof of Lemma \ref{lem3-3}}\\[.5ex]
(i) \ Consider the real-valued function $F$ on ${\cal X}^p\times \Theta$ defined by
\begin{equation}
F(z_1,\ldots,z_p;\theta)\,=\,\Bigl(\det\bigl[f_\theta(z_1),\ldots,f_\theta(z_p)\bigr]\Bigr)^2.
\label{eq3-5}
\end{equation}
Clearly, $F$ is continuous and hence, by compactness of ${\cal X}^p\times\Theta$, uniformly continuous.
Let
\begin{equation}
c_0\,=\,\min_{\theta\in\Theta}\max_{z_1,\ldots,z_p\in{\cal X}}F(z_1,\ldots,z_p;\theta).
\label{eq3-6}
\end{equation}
In fact, continuity of $\theta\longmapsto F(z_1,\ldots,z_p;\theta)$ for every fixed $(z_1,\ldots,z_p)\in{\cal X}^p$ 
implies lower semi-continuity of the function 
$\theta\longmapsto\max_{z_1,\ldots,z_p\in{\cal X}}F(z_1,\ldots,z_p;\theta)$, and hence this function attains its minimum on 
$\Theta$. The function is strictly positive on $\Theta$ (by the basic assumption (b4), (i)), hence its
minimum value  is positive, i.e., $c_0>0$. By the uniform continuity of $F$  
there exists a $\Delta_0>0$ such that
\begin{eqnarray}  
&&\max_{\theta\in\Theta}\big|F(z_1,\ldots,z_p;\theta)-F(z_1',\ldots,z_p';\theta)\big|\,<c_0
\nonumber\\
&&\mbox{for all $(z_1,\ldots,z_p),\ (z_1',\ldots,z_p')\in{\cal X}^p$ with ${\rm d}_{{\cal X}}(z_j,z_j')<\Delta_0$,
$1\le j\le p$.}\label{eq3-7}
\end{eqnarray}
 For any $k\ge1$ and $\ell,m\in\{1,\ldots,p\}$, $\ell<m$,  consider the particular points
$(x_{n_k+1},\ldots,x_{n_k+p})$ and $(z_1',\ldots,z_p')$ where the latter is given by
\[
z_j'=x_{n_k+j}\ \mbox{for $j\not=\ell$, and }\ z_\ell'=x_{n_k+m},
\]
and consider the parameter point $\theta=\theta_k$. Since the matrix 
$\bigl[f_{\theta_k}(z_1'),\ldots,f_{\theta_k}(z_p')\bigr]$ has two identical columns 
(the $\ell$-th and the $m$-th columns) one has $F(z_1'\ldots,z_p';\theta_k)=0$.
By (\ref{eq2-1}) $F(x_{n_k+1},\ldots,x_{n_k+p};\theta_k)=\max_{z_1,\ldots,z_p\in{\cal X}}F(z_1,\ldots,z_p;\theta_k)$,
and hence by (\ref{eq3-6})\\ 
$F(x_{n_k+1},\ldots,x_{n_k+p};\theta_k)\ge c_0$. 
Together with  $F(z_1',\ldots,z_p';\theta_k)=0$ and ${\rm d}_{{\cal X}}(x_{n_k+j},z_j')=0$ for $j\not=\ell$,
and ${\rm d}_{{\cal X}}(x_{n_k+\ell},z_\ell')={\rm d}_{{\cal X}}(x_{n_k+\ell},x_{n_k+m})$, ones gets from 
(\ref{eq3-7}) that ${\rm d}_{{\cal X}}(x_{n_k+\ell},x_{n_k+m})\ge\Delta_0$, which proves part (i) of the lemma.\\[.5ex]
(ii) \  The function $F$ from the proof of part (i) may be written as
\[
F(z_1,\ldots,z_p;\theta)\,=\,p^p\,\det\Bigl(p^{-1}\bigl[f_{\theta}(z_1),\ldots,f_{\theta}(z_p)\bigr]\,
\bigl[f_{\theta}(z_1),\ldots,f_{\theta}(z_p)\bigr]^\trp\Bigr),
\]
and the matrix on the r.h.s.~under the determinant is equal to the information matrix\\ 
$M\bigl(\xi[z_1,\ldots,z_p],\theta\bigr)$
where $\xi[z_1,\ldots,z_p]:=\frac{1}{p}\sum_{j=1}^p\delta[z_j]$. So (\ref{eq3-6}) rewrites as
\begin{equation}
p^{-p}c_0\,=\,\min_{\theta\in\Theta}\max_{z_1,\ldots,z_p\in{\cal X}}\det M\bigl(\xi[z_1,\ldots,z_p],\theta\bigr).
\label{eq3-6a}
\end{equation}
For $k\in\mathbb{N}$ the design $\eta_k=\frac{1}{p}\sum_{j=1}^p\delta[x_{n_k+j}]$ has the property  
\[
\det M(\eta_k,\theta_k) =\max_{z_1,\ldots,z_p\in{\cal X}} \det M\bigl(\xi[z_1,\ldots,z_p],\theta_k\bigr).
\]
Hence, by (\ref{eq3-6a}), 
\begin{equation}
\det M(\eta_k,\theta_k)\,\ge\,p^{-p}\,c_0\ \mbox{ for all $k\in\mathbb{N}$.}
\label{eq3-8}
\end{equation}
Using the positive finite constant $\gamma=\max_{(x,\theta)\in{\cal X}\times\Theta}\Vert f_\theta(x)\Vert$ it is easily seen that 
for all designs $\xi\in\Xi$, all $\theta\in\Theta$, and all vectors $a\in\mathbb{R}^p$ one has 
$a^\trp M(\xi,\theta)\,a\,\le\,\gamma^2\,a^\trp a$, where the identity
\begin{equation}
a^\trp M(\xi,\theta)\,a\,=\,\sum_{x\in{\rm\scriptsize supp}(\xi)}\xi(x)\,\bigl(a^\trp f_\theta(x)\bigr)^2
\label{eq3-8a}
\end{equation}
is useful. So $M(\xi,\theta)\le\gamma^2I_p$ in the Loewner semi-ordering, where
$I_p$ denotes the $p\times p$ identity matrix. Hence all the eigenvalues of $M(\xi,\theta)$ are less than or equal to $\gamma^2$.
Together with  (\ref{eq3-8}) it follows that
\begin{equation}
\lambda_{\rm\scriptsize min}\bigl(M(\eta_k,\theta_k)\bigr)\,\ge p^{-p}\,c_0\,\gamma^{-2(p-1)}\ \mbox{ for all $k\in\mathbb{N}$.}
\label{eq3-9}
\end{equation}
For any vector $a\in\mathbb{R}^p$ with $\Vert a\Vert=1$ one has
 $a^\trp M(\eta_k,\theta_k)$\,$a\,\ge \lambda_{\rm\scriptsize min}\bigl(M(\eta_k,\theta_k)\bigr)$, and hence
\[
a^\trp M(\eta_k,\theta_k)\,a\,\ge p^{-p}\,c_0\,\gamma^{-2(p-1)},
\]
and observing (\ref{eq3-8a}) this yields
\begin{equation}
\max_{1\le j\le p}\bigl(a^\trp f_{\theta_k}(x_{n_k+j})\bigr)^2\,\ge\,p^{-p}\,c_0\,\gamma^{-2(p-1)}.
\label{eq3-10}
\end{equation}
By (GLM) $f_\theta(x)=\psi(x,\theta)\,f(x)$ for all $(x,\theta)\in{\cal X}\times\Theta$ where, in particular, 
$\psi$ is a continuous positive function. So, $\psi_{\rm\scriptsize max}:=\max_{(x,\theta)\in{\cal X}\times\Theta}\psi(x,\theta)$
is a positive finite constant, and hence by (\ref{eq3-10}), defining \ 
$\varepsilon_0:= p^{-p/2}\,c_0^{1/2}\,\gamma^{-(p-1)}\psi_{\rm\scriptsize max}^{-1}$ one gets from (\ref{eq3-10})
\[
\max_{1\le j\le p}\bigl(a^\trp f(x_{n_k+j})\bigr)^2\,\ge\,\varepsilon_0^2\ 
\mbox{for all $a\in\mathbb{R}^p$, $\Vert a\Vert=1$, and all $k\in\mathbb{N}$.}
\]
Clearly, this is the same as 
\begin{equation}
\eta_k\Bigl(\bigl\{x\in{\cal X}\,:\,|a^\trp f(x)|\ge\varepsilon_0\bigr\}\Bigr)\,\ge\,1/p
\ \mbox{ for all $a\in\mathbb{R}^p$, $\Vert a\Vert=1$, and all $k\in\mathbb{N}$.}
\label{eq3-11}
\end{equation}
From (\ref{eq3-11})  and \ $\xi_k=\frac{n_1}{n_k}\xi_1+\frac{p}{n_k}\sum_{j=1}^{k-1}\eta_j$ for all $k\ge2$ 
according to (\ref{eq2-1}), one gets for all $a\in\mathbb{R}^p$, $\Vert a\Vert=1$, and all $k\ge1$,
\[
\xi_k\Bigl(\bigl\{x\in{\cal X}\,:\,|a^\trp f(x)|\ge\varepsilon_0\bigr\}\Bigr)\,\ge\,
\frac{p}{n_k}\frac{k-1}{p}\,=\,\frac{k-1}{n_k},
\]
where the inequality is trivial for $k=1$.
\eop

\subsection{Proofs of the theorems}
{\bf Proof of Theorem \ref{theo4-1}}\\[.5ex] 
We proceed basically as in the proof of Theorem 1 in \cite{FF-NG-RS-19}, with apppropriate 
modifications. Consider the random variables
\[
S_{n_k}(\theta):=\sum_{i=1}^{n_k}\bigl(Y_i-\mu(X_i,\theta)\bigr)^2\quad \mbox{and}\quad 
D_{n_k}(\theta,\overline{\theta})\,:=\,\sum_{i=1}^{n_k}\bigl(\mu(X_i,\theta)-\mu(X_i,\overline{\theta})\bigr)^2.
\]
for all $k\in\mathbb{N}$ and $\theta\in\Theta$. 
The least squares estimator $\widehat{\theta}_k^{(LS)}$ minimizes $S_{n_k}(\theta)$ over
$\theta\in\Theta$. 
For $\varepsilon>0$ we denote $C(\overline{\theta},\varepsilon):=
\bigl\{\theta\in\Theta\,:\,{\rm d}_\Theta(\theta,\overline{\theta})\ge\varepsilon\bigr\}$,
where ${\rm d}_\Theta$ denotes the distance function in $\Theta$.
The proof is divided into three steps.\\
\underbar{Step 1.} Show that for all $\varepsilon>0$ with 
$C(\overline{\theta},\varepsilon)\not=\emptyset$,
\[
\Big|\,\frac{1}{k}\Bigl(\inf_{\theta\in C(\overline{\theta},\varepsilon)} S_{n_k}(\theta)\,-S_{n_k}(\overline{\theta})\Bigr)\,
-\,\frac{1}{k}\inf_{\theta\in C(\overline{\theta},\varepsilon)} D_{n_k}(\theta,\overline{\theta})\,\Big|
\ \asto\,0\ \mbox{ (as $k\to\infty$).}
\]
\underbar{Step 2.} Show that for all $\varepsilon>0$ with 
$C(\overline{\theta},\varepsilon)\not=\emptyset$,
\[
\liminf_{k\to\infty}\Bigl(\frac{1}{k}\,
\inf_{\theta\in C(\overline{\theta},\varepsilon)}D_{n_k}(\theta,\overline{\theta})\Bigr)\ >\,0\ \ \mbox{a.s.}
\]
\underbar{Step 3.} Conclude from the results of Step 1 and Step 2 that for all $\varepsilon>0$ with 
$C(\overline{\theta},\varepsilon)\not=\emptyset$, 
\begin{equation}
\inf_{\theta\in C(\overline{\theta},\varepsilon)}S_{n_k}(\theta)\,-\,S_{n_k}(\overline{\theta})\ \asto\,\infty
\ \mbox{ (as $k\to\infty$).}
\label{eq4-1}
\end{equation}
Then, by (\ref{eq4-1}) and by Lemma 1 of Wu \cite{Wu} one gets \ 
$\widehat{\theta}_k^{(\rm\scriptsize LS)}\asto\overline{\theta}$.\\
\underbar{Ad Step 1.} 
As in \cite{FF-NG-RS-19} one gets, for all $k\in\mathbb{N}$,
\[
S_{n_k}(\theta)-S_{n_k}(\overline{\theta})\,=\,D_{n_k}(\theta,\overline{\theta})\,+\,2W_{n_k}(\theta,\overline{\theta}),\ \mbox{where }
W_{n_k}(\theta,\overline{\theta})\,:=\,\sum_{i=1}^{n_k}
\bigl(\mu(X_i,\overline{\theta})-\mu(X_i,\theta)\bigr)\,e_i.
\]
So, for all $k\ge1$,
\begin{equation}
\Big|\,\frac{1}{k}\Bigl(\inf_{\theta\in C(\overline{\theta},\varepsilon)} S_{n_k}(\theta)\,-S_{n_k}(\overline{\theta})\Bigr)\,
-\,\frac{1}{k}\inf_{\theta\in C(\overline{\theta},\varepsilon)} D_{n_k}(\theta,\overline{\theta})\,\Big|
\le\ \frac{2}{k}\,\sup_{\theta\in\Theta}\big|W_{n_k}(\theta,\overline{\theta})\bigr|.
\label{eq4-2}
\end{equation}
Introducing the function $h(x,\theta):=\mu(x,\overline{\theta})-\mu(x,\theta)$, $(x,\theta)\in{\cal X}\times\Theta$, 
we may write
\[
W_{n_k}(\theta,\overline{\theta})\,=\,
\sum_{j=1}^k\sum_{\ell=1}^p h(X_{n_{j-1}+\ell},\theta)\,e_{n_{j-1}+\ell},
\]
where for simplicity of presentation we assume that $n_1=p$,
and hence
\begin{equation}
\big|W_{n_k}(\theta,\overline{\theta})\big|\,\le\,
\sum_{\ell=1}^p\Big|\sum_{j=1}^k h(X_{n_{j-1}+\ell},\theta)\,e_{n_{j-1}+\ell}\Big|.
\label{eq4-2a}
\end{equation}
For each fixed $\ell\in\{1,\ldots,p\}$ the sequences of  random variables $e_{n_{j-1}+\ell}$, $j\in\mathbb{N}$,
and $X_{n_{j-1}+\ell}$, $j\in\mathbb{N}$, satisfy the assumptions of Lemma A.1 in \cite{FF-NG-RS-18}, and part (iii)
of that lemma yields
\begin{equation}
\frac{1}{k}\sup_{\theta\in\Theta}\Big|\sum_{j=1}^k h(X_{n_{j-1}+\ell},\theta)\,e_{n_{j-1}+\ell}\Big|\ \asto0
\ \mbox{ for each $\ell=1,\ldots,p$.}
\label{eq4-3}
\end{equation}
By (\ref{eq4-2}), (\ref{eq4-2a}) and (\ref{eq4-3}) the result of Step 1 follows.\\
\underbar{Ad Step 2 in the case that condition (SI) holds.}\\  
Consider any path $x_i,y_i$, $i\in\mathbb{N}$, and $\theta_k$, $k\in\mathbb{N}$
of the  sequences $X_i,Y_i$, $i\in\mathbb{N}$, and $\widehat{\theta}_k$, $k\in\mathbb{N}$. 
Choose $\Delta_0>0$  according to Lemma \ref{lem3-3} (i). 
Consider the subset of ${\cal X}^p$ given by
\[
D\,=\,\bigl\{(z_1,\ldots,z_p)\in{\cal X}^p\,:\,{\rm d}_{{\cal X}}(z_\ell,z_m)\ge\Delta_0,\ 1\le\ell<m\le p\bigr\},
\]
which is compact. Let $\varepsilon>0$ be given such that $C(\overline{\theta},\varepsilon)\not=\emptyset$.
Consider the (continuous) function on $D\times C(\overline{\theta},\varepsilon)$ given by
\[
(z_1,\ldots,z_p,\theta)\longmapsto
\sum_{\ell=1}^p\bigl(\mu(z_\ell,\theta)-\mu(z_\ell,\overline{\theta})\bigr)^2.
\]
By condition (SI) this function is strictly positive on $D\times C(\overline{\theta},\varepsilon)$,
and by compactness of this set the infimum $c(\varepsilon)$ of this function over $D\times C(\overline{\theta},\varepsilon)$
is attained and hence $c(\varepsilon)>0$. It follows that
\[
\sum_{\ell=1}^p\bigl(\mu(x_{n_k+\ell},\theta)-\mu(x_{n_k+\ell},\overline{\theta})\bigr)^2\,\ge\,c(\varepsilon)\ \ 
\mbox{for all $k\ge1$ and all $\theta\in  C(\overline{\theta},\varepsilon)$.}
\]
From this one gets for all $k\ge2$ and all $\theta\in C(\overline{\theta},\varepsilon)$, 
\begin{eqnarray*}
&&\frac{1}{k}D_{n_k}(\theta,\overline{\theta}) \ge  
\frac{1}{k}\sum_{i=n_1+1}^{n_k}\bigl(\mu(x_i,\theta)-\mu(x_i,\overline{\theta})\bigr)^2
=\frac{1}{k}\sum_{j=1}^{k-1}\sum_{\ell=1}^p \bigl(\mu(x_{n_j+\ell},\theta)-\mu(x_{n_j+\ell},\overline{\theta})\bigr)^2\\
&&\,\ge  \frac{(k-1)c(\varepsilon)}{k}\,\ge\,c(\varepsilon)/2.
\end{eqnarray*}
Hence the result of Step 2 follows.\\
\underbar{Ad Step 2 in the case that condition (${\rm GLM}^*$) holds.}\\ 
Again, consider any path $x_i,y_i$, $i\in\mathbb{N}$, and $\theta_k$, $k\in\mathbb{N}$
of the  sequences $X_i,Y_i$, $i\in\mathbb{N}$, and $\widehat{\theta}_k$, $k\in\mathbb{N}$. 
Choose a compact subinterval $J\subseteq I$ such that $f^\trp(x)\,\theta\in J$ for all $(x,\theta)\in\Theta$.
Then $b:=\min_{u\in J}G'(u)$ is positive and by the mean value theorem $|G(u)-G(v)|\ge b|u-v|$ for all
$u,v\in J$. Hence for all $i\in\mathbb{N}$ and $\theta\in\Theta$,
\[
\big|\mu(x_i,\theta)-\mu(x_i,\overline{\theta})\big|=
\big|G\bigl(f^\trp(x_i)\,\theta\bigr)-G\bigl(f^\trp(x_i)\,\overline{\theta}\bigr)\big|\,\ge\,
b\big|f^\trp(x_i)\,(\theta-\overline{\theta})\big|.
\]
So, for all $k\ge1$ and $\theta\in C(\overline{\theta},\varepsilon)$, denoting 
$a_\theta=(\theta-\overline{\theta})\big/\Vert\theta-\overline{\theta}\Vert$,
\begin{eqnarray}
&&D_{n_k}(\theta,\overline{\theta})\,=\,\sum_{i=1}^{n_k}\bigl(\mu(x_i,\theta)-\mu(x_i,\overline{\theta})\bigr)^2
\nonumber\\
&&\ge\,b^2\varepsilon^2\sum_{i=1}^{n_k}\bigl(f^\trp(x_i)\,a_\theta\bigr)^2\,=\,b^2\varepsilon^2n_k\int_{{\cal X}}
\bigl(f^\trp(x)\,a_\theta\bigr)^2\,{\rm d}\xi_k(x).\label{eq4-4}
\end{eqnarray}
Choose $\varepsilon_0>0$ according to Lemma \ref{lem3-3} (ii).  
Then for all $k\ge1$ and $\theta\in C(\overline{\theta},\varepsilon)$,
\begin{eqnarray*}
&&\int_{{\cal X}}\bigl(f^\trp(x)\,a_\theta\bigr)^2\,{\rm d}\xi_k(x)\,\ge
\int_{\bigl\{x\in{\cal X}\,:\,|f^\trp(x)\,a_\theta|\ge\varepsilon_0\big\}}\bigl(f^\trp(x)\,a_\theta\bigr)^2\,{\rm d}\xi_k(x)\\
&&\ge\,\varepsilon_0^2\,\xi_k\Bigl(\bigl\{x\in{\cal X}\,:\,|f^\trp(x)\,a_\theta|\ge\varepsilon_0\big\}\Bigr)\,\ge\,
\varepsilon_0^2(k-1)/n_k.
\end{eqnarray*}
Together with (\ref{eq4-4}) this yields
\begin{equation}
\frac{1}{k}\,\inf_{\theta\in C(\overline{\theta},\varepsilon)}D_{n_k}(\theta,\overline{\theta})\,\ge\,
\frac{1}{k}b^2\varepsilon^2n_k\varepsilon_0^2\frac{k-1}{n_k}.\label{eq4-5}
\end{equation}
For all $k\ge2$ the r.h.s. of (\ref{eq4-5}) is greater than or equal to 
$b^2\varepsilon^2\varepsilon_0^2/2$, and the result of Step 2 follows.\\
\underbar{Ad Step 3.}  Obviously, this follows from the results of steps 1 and 2. 
\eop

\noindent
{\bf Proof of Theorem \ref{theo4-2}}\\[.5ex] 
We will appropriately modify the arguments in the proof of Theorem 2 in \cite{FF-NG-RS-19}.
For simplicity of presentation, we assume $n_1=p$ for the starting design of the algorithm.
Choose a compact ball $\overline{B}$ centered at $\overline{\theta}$ such that $\overline{B}\subseteq{\rm int}(\Theta)$.
By the strong consistency of $\LSk$, $k\in\mathbb{N}$, there is a random variable $K$ with values in $\mathbb{N}\cup\{\infty\}$ such that $K<\infty$
a.s. and $\widehat{\theta}_k^{\rm\scriptsize (LS)}\in\overline{B}$ on $\{K\le k\}$ for all $k\in\mathbb{N}$. Along the lines 
in \cite{FF-NG-RS-19}, 
by equating the gradient w.r.t. $\theta$ of the sum of squares at $\theta=\LSk$ to zero, one gets for all $k\in\mathbb{N}$,
\begin{eqnarray}
\mbox{on $\{K\le k\}$\,: }\ \sum_{i=1}^{n_k}e_i\nabla\mu(X_i,\thb)\,=\, 
&&\sum_{i=1}^{n_k}\bigl[\mu(X_i,\LSk)-\mu(X_i,\thb)\bigr]\,\nabla\mu(X_i,\LSk)\nonumber\\
&& -\,\sum_{i=1}^{n_k}e_i\,\bigl[\nabla\mu(X_i,\LSk)-\nabla\mu(X_i,\thb)\bigr].
\label{eq4-6}
\end{eqnarray}
Concerning the asymtotic behavior of each of the three sums in (\ref{eq4-6}), 
we show the following. 
\begin{equation}
n_k^{-1/2}\sigma^{-1}(\thb)\,\Ms^{-1/2}(\thb)\sum_{i=1}^{n_k}e_i\nabla\mu(X_i,\thb) \dto {\rm N}(0,I_p) \ 
\mbox{ (as $k\to\infty$);} \label{eq4-7}
\end{equation}
\begin{eqnarray}
&&n_k^{-1/2}\sum_{i=1}^{n_k}\bigl[\mu(X_i,\LSk)-\mu(X_i,\thb)\bigr]\,\nabla\mu(X_i,\LSk)\nonumber\\
&&=\,\bigl[M(\xi_k,\LSk)\,+\,A_k\bigr]\,\bigl[n_k^{1/2}\bigl(\LSk-\thb\bigr)\bigr],\label{eq4-8}\\
&&\ \mbox{ with a sequence $A_k$, $k\in\mathbb{N}$, of random $p\times p$ matrices such that $A_k\asto0$;}\nonumber
\end{eqnarray}
\begin{eqnarray}
&&n_k^{-1/2}\sum_{i=1}^{n_k}e_i\bigl[\nabla\mu(X_i,\LSk)-\nabla\mu(X_i,\thb)\bigr]\,=\,
B_k\,\bigl[n_k^{1/2}\bigl(\LSk-\thb\bigr)\bigr],\label{eq4-9}\\
&&\ \mbox{ with a sequence $B_k$ of $p\times p$ random matrices such that $B_k\asto0$.}\nonumber
\end{eqnarray}
\underbar{Ad (\ref{eq4-7}).} 
\[
\sum_{i=1}^{n_k}e_i\nabla\mu(X_i,\thb)\,=\,
\sum_{j=1}^k\sum_{\ell=1}^pe_{n_{j-1}+\ell}\nabla\mu(X_{n_{j-1}+\ell},\thb)\,=\,\sum_{j=1}^k\overline{G}(\BM{X}_j)\,\BM{e}_j,
\]
where $\overline{G}$ denotes the $\mathbb{R}^{p\times p}$-valued function on ${\cal X}^p$ given by $\overline{G}(\BM{z})=
\bigl[\nabla\mu(z_1,\thb)\,,\,\ldots\,,\,\nabla\mu(z_p,\thb)\bigr]$ for all $\BM{z}=(z_1,\ldots,z_p)\in{\cal X}^p$.
Let any $v\in\mathbb{R}^p$ with $\Vert v\Vert=1$ be given.
Then,
\begin{eqnarray}
&&n_k^{-1/2}\sigma^{-1}(\thb)\,v^\trp \Ms^{-1/2}(\thb)\sum_{i=1}^{n_k}e_i\nabla\mu(X_i,\thb)\,=\,
k^{-1/2}\sum_{j=1}^k\widetilde{e}_j,\label{eq4-10}\\
&&\mbox{ where }\ \widetilde{e}_j=\BM{Z}_j^\trp\BM{e}_j\ \mbox{ and }\ 
\BM{Z}_j\,=\,p^{-1/2}\sigma^{-1}(\thb)\,\overline{G}^\trp(\BM{X}_j)\,\Ms^{-1/2}(\thb)\,v.\label{eq4-10a}
\end{eqnarray}
Note that $n_k=kp$ has been used. The sequence of $p$-dimensional random variables 
$\BM{Z}_j$, $j\in\mathbb{N}$, is uniformly bounded, that is, $\sup_{j\in\mathbb{N}}\Vert \BM{Z}_j\Vert\,\le c$ 
for some finite constant $c$, and $\BM{Z}_j$ is ${\cal F}_{j-1}$-measurable. From this it is easily seen that  
$\sum_{j=1}^k\widetilde{e}_j$, $k\in\mathbb{N}$, together with ${\cal F}_k$, $k\in\mathbb{N}\cup\{0\}$, 
constitutes a martingale. According to Corollary 3.1 of Hall and Heyde \cite{Hall-Heyde}, the distributional convergence 
$k^{-1/2}\sum_{j=1}^k\widetilde{e}_j\dto{\rm N}(0,1)$ is ensured by the following conditions $(\alpha)$ and $(\beta)$.\\[.5ex]
$(\alpha)$ \ $\displaystyle\frac{1}{k}\sum_{j=1}^k{\rm E}\bigl(\widetilde{e}_j^2\big|{\cal F}_{j-1}\bigr)\,\asto1$,\quad 
$(\beta)$ \ $\displaystyle\frac{1}{k}\sum_{j=1}^k
{\rm E}\Bigl(\widetilde{e}_j^2\Ifkt\bigl(|\widetilde{e}_j|>\varepsilon\sqrt{k}\bigr)\,\big|{\cal F}_{j-1}\Bigr)\,\asto0$.\\[.5ex]  
To verify $(\alpha)$, we write $\widetilde{e}_j^2=\BM{Z}_j^\trp\BM{e}_j\BM{e}_j^\trp\BM{Z}_j$, hence 
\begin{equation}
{\rm E}\bigl(\widetilde{e}_j^2\big|{\cal F}_{j-1}\bigr)\,=\,
\BM{Z}_j^\trp{\rm E}\bigl(\BM{e}_j\BM{e}_j^\trp\big|{\cal F}_{j-1}\bigr)\,\BM{Z}_j.\label{eq4-10b}
\end{equation}
By (AH) and the uniform boundedness of the sequence $\BM{Z}_j$, $j\in\mathbb{N}$,
\[
\max_{1\le j\le k}\Big|\BM{Z}_j^\trp{\rm E}\bigl(\BM{e}_j\BM{e}_j^\trp\big|{\cal F}_{j-1}\bigr)\,\BM{Z}_j\,-\,
\sigma^2(\thb)\,\BM{Z}_j^\trp\BM{Z}_j\Big|\,\asto0\ \mbox{ (as $k\to\infty$),}
\]
and hence
\begin{equation}
\Big|\frac{1}{k}\sum_{j=1}^k\BM{Z}_j^\trp{\rm E}\bigl(\BM{e}_j\BM{e}_j^\trp\big|{\cal F}_{j-1}\bigr)\,\BM{Z}_j\,-\,
\frac{\sigma^2(\thb)}{k}\sum_{j=1}^k\BM{Z}_j^\trp\BM{Z}_j\Big|\,\asto 0.\label{eq4-11}
\end{equation}
By (\ref{eq4-10a}),
\begin{equation}
\frac{\sigma^2(\thb)}{k}\sum_{j=1}^k\BM{Z}_j^\trp\BM{Z}_j\,=\,
v^\trp\Ms^{-1/2}\Bigl(\frac{1}{kp}\sum_{j=1}^k\overline{G}(\BM{X}_j)\,\overline{G}^\trp(\BM{X}_j)\Bigr)\,\Ms^{-1/2}(\thb)\,v.
\label{eq4-12}
\end{equation}
Since  $\overline{G}(\BM{X}_j)\,\overline{G}^\trp(\BM{X}_j)=
\sum_{\ell=1}^p\nabla\mu(X_{n_{k-1}+\ell},\thb)\,\nabla^\trp\mu(X_{n_{k-1}+\ell},\thb)$,
and by (b5),
\begin{equation}
\frac{1}{kp}\sum_{j=1}^k\overline{G}(\BM{X}_j)\,\overline{G}^\trp(\BM{X}_j)\,=\,
M(\xi_k,\thb).\label{eq4-13}
\end{equation}
By Corollary \ref{cor3-1}, $M(\xi_k,\thb)\asto \Ms(\thb)$ and hence, together with (\ref{eq4-13}), (\ref{eq4-12}),
(\ref{eq4-11}), and (\ref{eq4-10b}), condition $(\alpha)$ follows.\\
To verify $(\beta)$, we observe that 
$\widetilde{e}_j^2=\bigl(\BM{Z}_j^\trp\BM{e}_j\bigr)^2\le c^2\Vert\BM{e}_j\Vert^2$ and hence 
\[
{\rm E}\Bigl(\widetilde{e}_j^2\Ifkt\bigl(|\widetilde{e}_j|>\varepsilon\sqrt{k}\bigr)\,\big|{\cal F}_{j-1}\Bigr)\,
\le\, c^2\,
{\rm E}\Bigl(\Vert \BM{e}_j\Vert^2\Ifkt\bigl(\Vert\BM{e}_j\Vert>(\varepsilon/c)\sqrt{k}\bigr)\,\big|{\cal F}_{j-1}\Bigr).
\]
So $(\beta)$ follows from (L). We have thus shown that $k^{-1/2}\sum_{j=1}^k\widetilde{e}_j\dto{\rm N}(0,1)$, and together with (\ref{eq4-10})
and the Cram\'{e}r-Wold device, the convergence (\ref{eq4-7}) follows.\\[.5ex]
\underbar{Ad (\ref{eq4-8}).} \  As in \cite{FF-NG-RS-19}, one obtains
\begin{eqnarray*}
&&n_k^{-1/2}\sum_{i=1}^{n_k}\bigl[\mu(X_i,\LSk)-\mu(X_i,\thb)\bigr]\,\nabla\mu(X_i,\LSk)\,=\,
\bigl[M(\xi_k,\LSk)+A_k\bigr]\,\bigl[n_k^{1/2}(\LSk-\thb)\bigr],\\
&&\mbox{ where }\ A_k\,=\,\frac{1}{n_k}\sum_{i=1}^{n_k}\nabla\mu(X_i,\LSk)\,\bigl[
\nabla\mu(X_i,\widetilde{\theta}_{i,k})-\nabla\mu(X_i,\LSk)\bigr]^\trp,
\end{eqnarray*}
and where $\widetilde{\theta}_{i,k}$, $1\le i\le n_k$, are appropriate random points on the line segment joining $\LSk$ and $\thb$.
Along the lines in \cite{FF-NG-RS-19}, p.~11, one concludes $A_k\asto 0$.\\[.5ex]
\underbar{Ad (\ref{eq4-9}).} \  
Let any $v\in\mathbb{R}^p$ be given. As in \cite{FF-NG-RS-19} one calculates
\begin{eqnarray}
&&v^\trp\Bigl(n_k^{-1/2}\sum_{i=1}^{n_k}e_i\bigl[\nabla\mu(X_i,\LSk)-\nabla\mu(X_i,\thb)\bigr]\Bigr)\,=\,
b_k^\trp(v)\,\bigl[n_k^{1/2}(\LSk-\thb)\bigr], \label{eq4-14}\\
&&\mbox{ where }\ b_k(v)\,=\,\frac{1}{n_k}\sum_{i=1}^{n_k}e_i\nabla^2\mu(X_i,\widetilde{\theta}_{i,k})\,v,
\nonumber
\end{eqnarray}
with appropriate random points $\widetilde{\theta}_{i,k}$, $1\le i\le n_k$, on the line segment joining $\LSk$ and $\thb$,
and with the Hessians $\nabla^2\mu(x,\theta)$, $(x,\theta)\in{\cal X}\times{\rm int}(\Theta)$,  according to assumption (b5). 
We decompose
\begin{eqnarray*}
&&b_k(v)\,=\,b_k^{(1)}(v)\,+\,b_k^{(2)}(v),\\
&&\mbox{ where }\ b_k^{(1)}(v)\,=\,\frac{1}{n_k}\sum_{i=1}^{n_k}e_i\nabla^2\mu(X_i,\thb)\,v\\
&&\mbox{ and }\ \ b_k^{(2)}(v)\,=\,\frac{1}{n_k}\sum_{i=1}^{n_k}e_i
\bigl[\nabla^2\mu(X_i,\widetilde{\theta}_{i,k})\,v\,-\,\nabla^2\mu(X_i,\thb)\,v\bigr].
\end{eqnarray*}
Consider $b_k^{(1)}(v)$, which can be written as
\begin{equation}
b_k^{(1)}(v)\,=\,\frac{1}{p}\sum_{\ell=1}^p\Bigl(\frac{1}{k}\sum_{j=1}^ke_{n_{j-1}+\ell}
\nabla^2\mu(X_{n_{j-1}+\ell},\thb)\,v\Bigr).
\label{eq4-15}
\end{equation}
For fixed $\ell\in\{1,\ldots,p\}$, each component of the inner sum on the r.h.s. of (\ref{eq4-15}) satisfies the assumptions of 
Lemma A.1 in \cite{FF-NG-RS-18} and hence, by part (iii) of that lemma, 
converges almost surely to zero (as $k\to\infty$). By (\ref{eq4-15}) 
we conclude that $b_k^{(1)}(v)\asto0$. Consider $b_k^{(2)}(v)$. By the uniform continuity of $(x,\theta)\mapsto\nabla^2\mu(x,\theta)\,v$ on  
${\cal X}\times\overline{B}$ according to (b5), and by 
\[
\max_{1\le i\le n_k}\big\Vert \widetilde{\theta}_{i,k}-\thb\big\Vert\,\le\,\big\Vert \LSk-\thb\big\Vert\,\asto 0, 
\]
one gets
\[
\max_{1\le i\le n_k}\Big\Vert \nabla^2\mu(X_i,\widetilde{\theta}_{i,k})\,v\,-\,\nabla^2\mu(X_i,\thb)\,v\Big\Vert\,\asto 0.
\]
Since
\[
\bigl\Vert b_k^{(2)}(v)\big\Vert\,\le\,
\Bigl(\max_{1\le i\le n_k}\big\Vert \nabla^2\mu(X_i,\widetilde{\theta}_{i,k})\,v\,-\,\nabla^2\mu(X_i,\thb)\,v\big\Vert\Bigr)
\,\frac{1}{n_k}\sum_{i=1}^{n_k}|e_i|,
\]
the concergence $b_k^{(2)}(v)\asto 0$ will follow from $\limsup_{k\to\infty}\frac{1}{n_k}\sum_{i=1}^{n_k}|e_i|\,<\infty$ a.s.
In fact, 
\[
\frac{1}{n_k}\sum_{i=1}^{n_k}|e_i|\,=\,\frac{1}{p}\sum_{\ell=1}^p\Bigl(\frac{1}{k}\sum_{j=1}^k\big|e_{n_{j-1}+\ell}\big|\Bigr),
\]
and for each fixed $\ell\in\{1,\ldots,p\}$ an application of Lemma A.1, part (i) in \cite{FF-NG-RS-18} yields 
$\limsup_{k\to\infty}\frac{1}{k}\sum_{j=1}^k\big|e_{n_{j-1}+\ell}\big|\,<\infty$ a.s., and hence also
$\limsup_{k\to\infty}\frac{1}{n_k}\sum_{i=1}^{n_k}|e_i|\,<\infty$ a.s. We have thus shown that $b_k(v)\asto 0$.
Specializing $v$ to the elementary unit vectors $v^{(\ell)}$ of $\mathbb{R}^p$, $1\le \ell\le p$, and forming the
$p\times p$ matrix $B_k$ with rows $b_k^\trp(v^{(1)}),\ldots,b_k^\trp(v^{(\ell)})$ one has $B_k\asto 0$, and 
(\ref{eq4-9}) follows from (\ref{eq4-14}).\\
From (\ref{eq4-6}), (\ref{eq4-7}), (\ref{eq4-8}) and (\ref{eq4-9}),
\[
\sigma^{-1}(\thb)\,\Ms^{-1/2}(\thb)\,\bigl[M(\xi_k,\LSk) +A_k-B_k\bigr]\,\bigl[\sqrt{kp}\,(\LSk-\thb)\bigr]\,
\dto {\rm N}(0,I_p),
\]
and $A_k\asto0$, $B_k\asto0$. By Corollary \ref{cor3-1}, $M(\xi_k,\LSk)\asto\Ms(\thb)$. 
Using standard properties of convergence in distribution, the result follows.
\eop

\noindent
{\bf Proof of Theorem \ref{theo5-1}}\\[.5ex] 
The error variables in model (a1), (a2'), (a3') are given by
\begin{equation}
e_i\,=\,Y_i\,-\,G\bigl(f^\trp(X_i)\,\thb\bigr), \ \ i\in\mathbb{N},\label{eq5-3}
\end{equation}
and we consider the error vectors
\begin{equation}
\BM{e}_k\,=\,\bigl(e_{n_{k-1}+1},\ldots,e_{n_k}\bigr)^\trp,\ \ k\in\mathbb{N}.\label{eq5-4}
\end{equation}
From (a2') and from general properties of an exponential family one concludes in particular, that
the fourth conditional moments of the error vectors are bounded by some finite constant $C_4$,
\begin{equation}
{\rm E}\bigl(\Vert \BM{e}_k\Vert^4\big|{\cal F}_{k-1}\bigr)\,\le\,C_4\ \mbox{ a.s. for all $k\in\mathbb{N}$.}
\label{eq5-5}
\end{equation}
Along the lines of the proof of Theorem 3.3 in \cite{FF-NG-RS-18} one obtains, for all $\theta\in\Theta$
and all $k\in\mathbb{N}$,
\begin{equation}
L_{n_k}(\overline{\theta})-L_{n_k}(\theta)\,\ge\,\sum_{i=1}^{n_k}\bigl(\overline{\tau}_i-\tau_i(\theta)\bigr)e_i + 
{\textstyle\frac{1}{2}}\beta_0\beta_1^2\sum_{i=1}^n\bigl[f^\trp(X_i)\,\theta - f^\trp(X_i)\,\overline{\theta}\,\bigr]^2
\label{eq5-6}
\end{equation}
with some positive real constants $\beta_0$ and $\beta_1$. 
According to Wu \cite{Wu}, Lemma 1, for
strong consistency of $\widehat{\theta}_k^{\rm\scriptsize (ML)}$ it is sufficient to show that,
for every $\varepsilon>0$ such that the parameter subset 
$C(\overline{\theta},\varepsilon)=\bigl\{\theta\in\Theta\,:\,\Vert\theta-\overline{\theta}\Vert\ge\varepsilon\bigr\}$
is nonempty, one has
\[
\liminf_{k\to\infty}\Bigl(L_{n_k}(\overline{\theta})-\sup_{\theta\in C(\overline{\theta},\delta)} L_{n_k}(\theta)\Bigr)\,>0 \ \mbox{ a.s.}
\]
In fact, the $\liminf$ turns out to be equal to infinity almost surely, since we show that    
\begin{equation}
\liminf_{k\to\infty}\frac{1}{k}\Bigl(L_{n_k}(\overline{\theta})-\sup_{\theta\in C(\overline{\theta},\varepsilon)} L_{n_k}(\theta)\Bigr)\ >0
 \ \mbox{a.s.}\label{eq5-7}
\end{equation}
From (\ref{eq5-6}) one concludes
\begin{eqnarray}
&&\frac{1}{k}\Bigl(L_{n_k}(\overline{\theta})-\sup_{\theta\in C(\overline{\theta},\varepsilon)} L_{n_k}(\theta)\Bigr)
\nonumber\\
&&\ge\,-\frac{1}{k}\sup_{\theta\in\Theta}\Big|\sum_{i=1}^{n_k}\bigl(\overline{\tau}_i-\tau_i(\theta)\bigr)e_i\Big|
+ {\textstyle\frac{1}{2}}\beta_0\beta_1^2
\frac{1}{k}\inf_{\theta\in C(\overline{\theta},\varepsilon)}\sum_{i=1}^{n_k}\bigl[f^\trp(X_i)\,\theta - f^\trp(X_i)\,\overline{\theta}\,\bigr]^2.\label{eq5-8}
\end{eqnarray}                   
Introduce the function $h(x,\theta)\,:=\,(b')^{-1}\bigl(f^\trp(x)\,\overline{\theta}\bigr) - (b')^{-1}\bigl(f^\trp(x)\,\theta\bigr)$,
$(x,\theta)\in{\cal X}\times\Theta$. From the definition of $\overline{\tau}_i$ and $\tau_i(\theta)$ one has
$\overline{\tau}_i-\tau_i(\theta)=h(X_i,\theta)$ for all $i\in\mathbb{N}$ and $\theta\in\Theta$. 
For convenience, we now assume $n_1=p$. Then   
\begin{eqnarray}
&&\Big|\sum_{i=1}^{n_k}h(X_i,\theta)\,e_i\Big|\,=\,
\Big|\sum_{\ell=1}^p\sum_{j=1}^kh(X_{n_{j-1}+\ell},\theta)\,e_{n_{j-1}+\ell}\Big|\nonumber\\
&&\le\,\sum_{\ell=1}^p\Big|\sum_{j=1}^kh(X_{n_{j-1}+\ell},\theta)\,e_{n_{j-1}+\ell}\Big|.\label{eq5-9}
\end{eqnarray}
By (a1), (a2') and an application of Lemma A.1, part (iii) of \cite{FF-NG-RS-18}, one gets for each $\ell=1,\ldots,p$,
\[
\frac{1}{k}\sup_{\theta\in\Theta}\Big|\sum_{j=1}^kh(X_{n_{j-1}+\ell},\theta)\,e_{n_{j-1}+\ell}\Big|\,\asto 0\
\]
and hence by (\ref{eq5-9})
\begin{equation}
\frac{1}{k}\sup_{\theta\in\Theta}\Big|\sum_{i=1}^{n_k}h(X_i,\theta)\,e_i\Big|\,\le\,
\sum_{\ell=1}^p \Bigl(
\frac{1}{k}\sup_{\theta\in\Theta}\Big|\sum_{j=1}^kh(X_{n_{j-1}+\ell},\theta)\,e_{n_{j-1}+\ell}\Big|\,\Bigr)\ \asto 0.
\label{eq5-9a}
\end{equation}
In view of (\ref{eq5-8}) and (\ref{eq5-9a}), it remains to show that
\begin{equation}
\liminf_{k\to\infty}\Bigl(
\frac{1}{k}\inf_{\theta\in C(\overline{\theta},\varepsilon)}
\sum_{i=1}^{n_k}\bigl[f^\trp(X_i)\,\theta - f^\trp(X_i)\,\overline{\theta}\,\bigr]^2\Bigr)\ >\,0\ \mbox{ a.s.}
\label{eq5-10}
\end{equation}
To this end, we consider an arbitrary path of the adaptive process and, in particular, a path $x_i$, $i\in\mathbb{N}$, of the sequence
$X_i$, $i\in\mathbb{N}$. Since 
\[
\int_{{\cal X}}\bigl[f^\trp(x)\,\theta - f^\trp(x)\,\overline{\theta}\,\bigr]^2\,{\rm d}\xi_k(x)\,=\,
\frac{1}{n_k}\sum_{i=1}^{n_k}\bigl[f^\trp(x_i)\,\theta - f^\trp(x_i)\,\overline{\theta}\,\bigr]^2,
\]
(\ref{eq5-10}) will follow from  
\begin{equation}
\liminf_{k\to\infty}\Bigl(\inf_{\theta\in C(\overline{\theta},\varepsilon)}
 \int_{{\cal X}}\bigl[f^\trp(x)\,(\theta - \overline{\theta})\,\bigr]^2\,{\rm d}\xi_k(x)\,\Bigr)\ >0.
\label{eq5-11}
\end{equation}
In fact, (\ref{eq5-11}) can be seen as follows. By Lemma \ref{lem3-3}, part (ii), there is an $\varepsilon_0>0$ such that
for all $k\in\mathbb{N}$ and all $a\in\mathbb{R}^p$ with $\Vert a\Vert=1$ one has
\begin{equation}
\xi_k\Bigl(\bigl\{x\in{\cal X}\,:\,|f^\trp(x)\,a|\ge\varepsilon_0\bigr\}\Bigr)\,\ge\,(k-1)/n_k.
\label{eq5-12}
\end{equation}
In particular, for any $\theta\in C(\overline{\theta},\varepsilon)$ we take $a=a_\theta=(\theta-\thb)/\Vert\theta-\thb\Vert$,
and from (\ref{eq5-12}) together with $\Vert\theta-\thb\Vert\ge\varepsilon$ we get
\[
\xi_k\Bigl(\bigl\{x\in{\cal X}\,:\,|f^\trp(x)\,(\theta-\thb)|\ge\varepsilon_0\varepsilon\bigr\}\Bigr)\,\ge\,(k-1)/n_k
\]
for all $k\in\mathbb{N}$ and all $\theta\in C(\overline{\theta},\varepsilon)$. Hence, using the obvious inequality
\[
\int_{{\cal X}}\bigl[f^\trp(x)\,(\theta - \overline{\theta})\,\bigr]^2\,{\rm d}\xi_k(x)\,\ge\,
(\varepsilon_0\varepsilon)^2\,\xi_k\Bigl(\bigl\{x\in{\cal X}\,:\,|f^\trp(x)\,(\theta-\thb)|\ge\varepsilon_0\varepsilon\bigr\}\Bigr),
\]
we obtain
\[
\inf_{\theta\in C(\overline{\theta},\varepsilon)}
 \int_{{\cal X}}\bigl[f^\trp(x)\,(\theta - \overline{\theta})\,\bigr]^2\,{\rm d}\xi_k(x)\,\ge\,
(\varepsilon_0\varepsilon)^2(k-1)/n_k.
\]
Since $\lim_{k\to\infty}(k-1)/n_k\,=1/p$, it follows that the $\liminf$ in (\ref{eq5-11}) is greater than or equal to
$(\varepsilon_0\varepsilon)^2/p\,>0$.
\eop

\noindent
{\bf Proof of Theorem \ref{theo5-2}}\\[.5ex] 
We will appropriately modify the arguments in the proof of Theorem 3.3 in \cite{FF-NG-RS-18}.
For simplicity of presentation, we assume $n_1=p$ for the starting design of the algorithm.
Choose a compact ball $\overline{B}$ centered at $\overline{\theta}$ such that $\overline{B}\subseteq{\rm int}(\Theta)$.
By the strong consistency of $\MLk$ according to Theorem \ref{theo5-1}, 
there is a random variable $K$ with values in $\mathbb{N}\cup\{\infty\}$ such that $K<\infty$
a.s. and $\widehat{\theta}_k^{\rm\scriptsize (ML)}\in\overline{B}$ on $\{K\le k\}$ for all $k\in\mathbb{N}$. 
Along the lines of \cite{FF-NG-RS-18}, p.~719, one obtains for the gradients (w.r.t. $\theta$) of the log-likelihood,  
$S_{n_k}(\theta)=\nabla L_{n_k}(\theta)$, where $\theta\in\overline{B}$,
\begin{eqnarray}
&&S_{n_k}(\theta)\,=\,\sum_{i=1}^{n_k}\Bigl(Y_i-G\bigl(f^\trp(X_i)\,\theta\bigr)\Bigr)\,H\bigl(f^\trp(X_i)\,\theta\bigr)\,f(X_i),
\label{eq5-13}\\
&&\mbox{ where }\ H(u)\,=\,\frac{G'(u)}{b''\Bigl((b')^{-1}\bigl(G(u)\bigr)\Bigr)}\ \mbox{ for all $u\in I$,}
\label{eq5-13a}
\end{eqnarray}
and one concludes,
\begin{eqnarray}
&&\mbox{ a.s. on $\{K\le k\}$\,: } \quad \frac{1}{\sqrt{n_k}}\,\Ms^{-1/2}(\thb)\,S_{n_k}(\overline{\theta})\nonumber\\
&&\phantom{\mbox{ a.s. on xx\,: }}\,=\,\Ms^{-1/2}(\thb)\,\Bigl[M(\xi_k,\overline{\theta})+
\frac{1}{n_k}D_k - \frac{1}{n_k}B_k\Bigr]\,\bigl[\sqrt{n_k}\bigl(\widehat{\theta}_k^{\rm\scriptsize (ML)}-\overline{\theta}\bigr)\bigr],
 \label{eq5-14}\\[.5ex]
&&\mbox{where } \ 
D_k \,=\, \sum_{i=1}^{n_k}\bigl[\varphi^2\bigl(f^\trp(X_i)\,\widetilde{\theta}_{i,k}\bigr)-
\varphi^2\bigl(f^\trp(X_i)\,\overline{\theta}\bigr)\bigr]\,f(X_i)\,f^\trp(X_i)\label{eq5-14a}\\
&&\mbox{ \ and }\ 
B_k \,=\, \sum_{i=1}^{n_k}\Bigl(Y_i-G\bigl(f^\trp(X_i)\,\widetilde{\theta}_{i,n}\bigr)\Bigr)
\,H'\bigl(f^\trp(X_i)\,\widetilde{\theta}_{i,k}\bigr)\,f(X_i)\,f^\trp(X_i), \label{eq5-14b}
\end{eqnarray}
and where $\widetilde{\theta}_{i,k}$, $1\le i\le n_k$, are appropriate random points on the line segment joining $\LSk$ and 
$\thb$.
The asymptotics (as $k\to\infty$) of the left-hand side of (\ref{eq5-14}) and of the random matrices $D_k$ and $B_k$ 
will shown to be as follows.
\begin{eqnarray}
&& \frac{1}{\sqrt{n_k}}\,\Ms^{-1/2}(\thb)\,S_{n_k}(\overline{\theta})\,\dto\,{\rm N}(0,I_p);\label{eq5-15}\\
&& \frac{1}{n_k}D_k\asto 0\ \mbox{ and }\ \frac{1}{n_k}B_k\asto 0.\label{eq5-15ab}
\end{eqnarray}
\underline{Ad (\ref{eq5-15}):} \ 
According to the Cram\'{e}r-Wold device choose any $v\in\mathbb{R}^p$ with $\Vert v\Vert=1$. 
By (\ref{eq5-13}) with $\theta=\overline{\theta}$ and  $Y_i=G\bigl(f^\trp(X_i)\,\overline{\theta}\bigr)\,+e_i$
according to (\ref{eq5-3}), we get
\begin{eqnarray}
\frac{1}{\sqrt{n_k}}\,v^\trp \Ms^{-1/2}(\thb)\,S_n(\overline{\theta})&=&\,
\frac{1}{\sqrt{n_k}}\,v^\trp \Ms^{-1/2}(\thb)\sum_{i=1}^{n_k}
e_iH\bigl(f^\trp(X_i)\,\overline{\theta}\bigr)\,f(X_i)\nonumber\\
&=& \frac{1}{\sqrt{kp}}\,v^\trp \Ms^{-1/2}(\thb)\sum_{\ell=1}^p\sum_{j=1}^k\overline{G}(\BM{X}_j)\,\BM{e}_j,\label{eq5-16}
\end{eqnarray}
where \ $\overline{G}\,:\,{\cal X}^p\longrightarrow\mathbb{R}^{p\times p}$ is defined by
\[
\overline{G}(\BM{z})\,:=\,\Bigl[H\bigl(f^\trp(z_1)\,\thb\bigr)\,f(z_1)\,,\,\ldots\,,\,
 H\bigl(f^\trp(z_p)\,\thb\bigr)\,f(z_p)\Bigr]\ \mbox{ for }\BM{z}=(z_1,\ldots,z_p)\in{\cal X}^p,
\]
and the error vectors $\BM{e}_j$, $j\in\mathbb{N}$, are given by (\ref{eq5-4}).
Introducing the sequence $\BM{Z}_j$, $j\in\mathbb{N}$, of $\mathbb{R}^p$-valued random variables,
\begin{equation} 
\BM{Z}_j\,=\,p^{-1/2}\,\overline{G}^\trp(\BM{X}_j)\,\Ms^{-1/2}(\thb)\,v,\label{eq5-17}
\end{equation}
we can write
\begin{equation}
\frac{1}{\sqrt{n_k}}\,v^\trp \Ms^{-1/2}(\thb)\,S_n(\overline{\theta})\,=\,
\frac{1}{\sqrt{k}}\,\sum_{j=1}^k\BM{Z}_j^\trp\BM{e}_j.\label{eq5-18}
\end{equation}
Abbreviate $\widetilde{e}_j=\BM{Z}_j^\trp\BM{e}_j$. 
Since $\BM{Z}_j$ is ${\cal F}_{j-1}$-measurable for all $j\in\mathbb{N}$, and the sequence $\BM{Z}_j$
is uniformly bounded, that is, $\Vert\BM{Z}_j\Vert\le c$ for all $j\in\mathbb{N}$ for some finite constant $c$,
it follows that the sequence of partial sums $\sum_{j=1}^k\widetilde{e}_k$, $k\in\mathbb{N}$, is a martingale w.r.t.
the filtration ${\cal F}_k$, $k\in\mathbb{N}\cup\{0\}$. We will verify the following two conditions (1) and (2).\\[.5ex]  
(1) \ $\displaystyle\frac{1}{k}\sum_{j=1}^k{\rm E}\bigl(\widetilde{e}_j^2\big|\,{\cal F}_{j-1}\bigr)\,\asto 1$;
 \ \ (2) \ $\displaystyle\frac{1}{k}\sum_{j=1}^k
{\rm E}\Bigl(\widetilde{e}_j^2\Ifkt\bigl(|\widetilde{e}_j|>\sqrt{k}\,\varepsilon\bigr)\big|\,{\cal F}_{j-1}\Bigr)
\,\asto 0$ \ for all $\varepsilon>0$. \\[.5ex]
Then, by Corollary 3.1 (p.~58) of Hall and Heyde \cite{Hall-Heyde}, the convergence \  
$\frac{1}{\sqrt{k}}\,\sum_{j=1}^k\widetilde{e}_j\dto {\rm N}(0,1)$ and thus (\ref{eq5-15}) will follow.
To verify condition (1), inserting 
$\widetilde{e}_j^2=\bigl(\BM{Z}_j^\trp\BM{e}_j\bigr)^2=\BM{Z}_j^\trp\BM{e}_j\BM{e}_j^\trp\BM{Z}_j$, one gets
${\rm E}\bigl(\widetilde{e}_j^2\big|{\cal F}_{j-1}\bigr)=\BM{Z}_j^\trp{\rm E}\bigl(\BM{e}_j\BM{e}_j^\trp\big|{\cal F}_{j-1}\bigr)\,\BM{Z}_j$,
and hence
\[
 \frac{1}{k}\sum_{j=1}^k{\rm E}\bigl(\widetilde{e}_k^2\big|\,{\cal F}_{k-1}\bigr)\,=\,
\frac{1}{k}\sum_{j=1}^k\BM{Z}_j^\trp{\rm E}\bigl(\BM{e}_j\BM{e}_j^\trp\big|{\cal F}_{j-1}\bigr)\,\BM{Z}_j.
\]
According to (\ref{eq5-3}) and (a2'),
\[
{\rm E}\bigl(\BM{e}_j\BM{e}_j^\trp\big|{\cal F}_{j-1}\bigr)\,=\,\overline{V}(\BM{X}_j),
\]
where $\overline{V}(\BM{z})\,:=\,
{\rm diag}\Bigl(b''\bigl((b')^{-1}(f^\trp(z_\ell)\,\thb)\bigr)\,:\,1\le\ell\le p\Bigr)$ for $\BM{z}=(z_1,\ldots,z_p)\in{\cal X}^p$,
and where ${\rm diag}\bigl(a_\ell\,:\,1\le\ell\le p\bigr)$, for real numbers $a_1,\ldots, a_p$, 
stands for the diagonal $p\times p$ matrix with diagonal entries $a_1,\ldots, a_p$.
Inserting according to (\ref{eq5-17}), one gets
\[
\frac{1}{k}\sum_{j=1}^k{\rm E}\bigl(\widetilde{e}_j^2\big|\,{\cal F}_{j-1}\bigr)\,=\,
v^\trp\Ms^{-1/2}(\thb)\Bigl(\frac{1}{kp}\sum_{j=1}^k\overline{G}(\BM{X}_j)\,\overline{V}(\BM{X}_j)\,\overline{G}^\trp(\BM{X}_j)\Bigr)\,
\Ms^{-1/2}(\thb)\,v.
\]
By the definitions of $\overline{G}(\BM{z})$ and $\overline{V}(\BM{z})$, where $\BM{z}=(z_1,\ldots,z_p)\in{\cal X}^p$,
and by (\ref{eq5-13a}) and (a3'), one gets
\[
\overline{G}(\BM{z})\,\overline{V}(\BM{z})\,\overline{G}^\trp(\BM{z})\,=\,
\sum_{\ell=1}^p\varphi^2\bigl(f^\trp(z_\ell)\,\thb\bigr)\,f(z_\ell)\,f^\trp(z_\ell).
\]
It follows that
\[
\frac{1}{kp}\sum_{j=1}^k\overline{G}(\BM{X}_j)\,\overline{V}(\BM{X}_j)\,\overline{G}^\trp(\BM{X}_j)\,=\,
\frac{1}{n_k}\sum_{i=1}^{n_k}\varphi^2\bigl(f^\trp(X_i)\,\thb\bigr)\,f(X_i)\,f^\trp(X_i)\,=\,
M(\xi_k,\thb).
\]
By Corollary \ref{cor3-1}, $M(\xi_k,\thb)\asto\Ms(\thb)$ which entails $\Ms^{-1/2}(\thb)\,M(\xi_k,\thb)\,\Ms^{-1/2}(\thb)\asto I_p$, 
and hence 
\[
\frac{1}{k}\sum_{j=1}^k{\rm E}\bigl(\widetilde{e}_j^2\big|\,{\cal F}_{j-1}\bigr)\,=\,
v^\trp\Ms^{-1/2}(\thb)\,M(\xi_k,\thb)\,\Ms^{-1/2}(\thb)\,v\,\asto1. 
\]
To verify condition (2), recall (\ref{eq5-5}) showing boundedness of the fourth conditional moments 
of the error vectors $\BM{e}_j$, $j\in\mathbb{N}$, by  a finite constant $C_4$, and recall also the uniform 
boundedness of the random vectors $\BM{Z}_j$, $j\in\mathbb{N}$, by a finite constant $c$.
Using the inequalities\\ 
$\widetilde{e}_j^2\Ifkt(|\widetilde{e}_j|>\sqrt{k}\varepsilon)\le\frac{1}{\varepsilon^2k}\widetilde{e}_j^4$
\ and $\widetilde{e}_j^4=(\BM{Z}^\trp_j\BM{e}_j)^4\le \Vert\BM{Z}_j\Vert^4\Vert\BM{e}_j\Vert^4$,
one gets
\[
\frac{1}{k}\sum_{j=1}^k
{\rm E}\Bigl(\widetilde{e}_j^2\Ifkt\bigl(|\widetilde{e}_j|>\sqrt{k}\,\varepsilon\bigr)\big|\,{\cal F}_{j-1}\Bigr)
\,\le\,\frac{1}{\varepsilon^2k^2}\sum_{j=1}^k\Vert\BM{Z}_j\Vert^4{\rm E}\bigl(\Vert\BM{e}_j\Vert^4\big|{\cal F}_{j-1}\bigr)
\le \frac{c^4C_4}{\varepsilon^2k} \ \to 0 \ \mbox{ (as $k\to\infty$).}
\]
\underline{Ad (\ref{eq5-15ab}):} \  
The first convergence $D_k/n_k\asto0$ is shown as in \cite{FF-NG-RS-18}, pp.~720-721. For the second convergence the arguments
in \cite{FF-NG-RS-18}, p.~721,  are modified as follows. We split $B_k$,
\begin{eqnarray*}
B_k &=&\,B_k^{(1)} + B_k^{(2)} + B_k^{(3)},\ \mbox{ where}\\
B_k^{(1)} &=& \sum_{i=1}^{n_k}\bigl[G\bigl(f^\trp(X_i)\,\overline{\theta}\bigr) -G\bigl(f^\trp(X_i)\,\widetilde{\theta}_{i,k}\bigr)\bigr]\,
H'\bigl(f^\trp(X_i)\,\widetilde{\theta}_{i,k}\bigr)\,f(X_i)\,f^\trp(X_i),\\
B_k^{(2)} &=& \sum_{i=1}^{n_k} e_i\,H'\bigl(f^\trp(X_i)\,\overline{\theta}\bigr)\,f(X_i)\,f^\trp(X_i),\\
B_k^{(3)} &=& \sum_{i=1}^{n_k} e_i\,
\bigl[H'\bigl(f^\trp(X_i)\,\widetilde{\theta}_{i,k}\bigr)-H'\bigl(f^\trp(X_i)\,\overline{\theta}\bigr)\bigr]\,
f(X_i)\,f^\trp(X_i).
\end{eqnarray*}
The convergence $\frac{1}{n_k}B_k^{(1)}\asto 0$ is obtained as in \cite{FF-NG-RS-18}, p.~721. 
Concerning $B_k^{(2)}$, fix any pair $(r,s)$,
where $1\le r,s\le p$. Consider the $(r,s)$-th entry of $\frac{1}{n_k}B_k^{(2)}$, which can be written as
\begin{eqnarray*}
&& \frac{1}{p}\sum_{\ell=1}^p\Bigl(\frac{1}{k}\sum_{j=1}^kZ_{j,\ell}\,e_{n_{j-1}+\ell}\Bigr),\\
&&\mbox{where }\ Z_{j,\ell}\,:=\,H'\bigl(f^\trp(X_{n_{j-1}+\ell})\,\overline{\theta}\bigr)\,
f_r(X_{n_{j-1}+\ell})\,f_s(X_{n_{j-1}+\ell}),\ \ j\in\mathbb{N},\ 1\le\ell\le p. 
\end{eqnarray*}
For each fixed $\ell\in\{1,\ldots,p\}$ an application of Lemma A.1, part (ii), of \cite{FF-NG-RS-18} yields
$\frac{1}{k}\sum_{j=1}^kZ_{j,\ell}\,e_{n_{j-1}+\ell}\asto0$. Hence each entry of $\frac{1}{n_k}B_k^{(2)}$ converges to zero
almost surely,
that is, $\frac{1}{n_k}B_k^{(2)}\asto0$.  Concerning $B_k^{(3)}$, as in \cite{FF-NG-RS-18}, p.~721, it is easily seen that 
the absolute value of each entry of $\frac{1}{n_k}B_k^{(3)}$  is bounded above by
\[
\gamma_0^2\,\Bigl(\max_{1\le i\le n_k}
\big|H'\bigl(f^\trp(X_i)\,\widetilde{\theta}_{i,k}\bigr)-H'\bigl(f^\trp(X_i)\,\overline{\theta}\bigr)\big|\Bigr)\,
\frac{1}{n_k}\sum_{i=1}^{n_k}|e_i|,
\]
and 
\[
\max_{1\le i\le n_k}
\big|H'\bigl(f^\trp(X_i)\,\widetilde{\theta}_{i,n}\bigr)-H'\bigl(f^\trp(X_i)\,\overline{\theta}\bigr)\big|\,\asto0. 
\]
Writing
\[
\frac{1}{n_k}\sum_{i=1}^{n_k}|e_i|\,=\,\frac{1}{p}\sum_{\ell=1}^p\Bigl(\frac{1}{k}\sum_{j=1}^k|e_{n_{j-1}+\ell}|\Bigr),
\]
an application of Lemma A.1, part (i), of \cite{FF-NG-RS-18} yields for each fixed $\ell\in\{1,\ldots,p\}$ that
\[
\limsup_{k\to\infty}\frac{1}{k}\sum_{j=1}^k|e_{n_{j-1}+\ell}|\,<\infty\ \mbox{ a.s.}
\]
Hence, $\limsup_{k\to\infty}\frac{1}{n_k}\sum_{i=1}^{n_k}|e_i|\,<\infty$ a.s., and thus each entry of $\frac{1}{n_k}B_k^{(3)}$
converges to zero almost surely, that is, $\frac{1}{n_k}B_k^{(3)}\asto0$.\\
From (\ref{eq5-14}), (\ref{eq5-15}), (\ref{eq5-15ab}), together with $M(\xi_k,\thb)\asto\Ms(\thb)$ and $\Ifkt(K\le k)\asto1$,
one concludes that 
\[
\Ms^{1/2}(\thb)\,\bigl[\sqrt{kp}\,(\MLk-\thb)\bigr]\dto {\rm N}(0,I_p),\ \mbox{ or equivalently, }\ 
\sqrt{kp}\,(\MLk-\thb)\dto {\rm N}\bigl(0,\Ms^{-1}(\thb)\bigr).
\]
\eop
\end{appendix}


\end{document}